\documentclass[12pt]{article}
\usepackage[english]{babel}
\usepackage{amsthm,amsfonts, amsbsy, amssymb,amsmath,graphicx}
\usepackage{graphics}
\usepackage{pifont}
\usepackage{amsthm}
\usepackage{amsfonts}
\usepackage{amssymb}
\usepackage[mathscr]{euscript}
\usepackage[all]{xy}

\def\thtext#1{
\catcode`@=11 \gdef\@thmcountersep{. #1}
\catcode`@=12}

\renewcommand{\baselinestretch}{1.2}

\newtheorem{theorem}{Theorem}[section]

\newcounter{imn1}[section]
\renewcommand{\thtext}
{\thesection.\arabic{imn1}}

\newcounter{imn2}[section]
\renewcommand{\thtext}
{\thesection.\arabic{imn2}}

\newcounter{imn3}[section]
\renewcommand{\thtext}
{\thesection.\arabic{imn3}}

 \newenvironment{definition}{\trivlist \item[\hskip\labelsep{\bf Definition}]
 \refstepcounter{imn1}{\bf\thesection.\arabic{imn1}.}}%
 {\endtrivlist}

 {\endtrivlist}

 \newenvironment{exa}{\trivlist \item[\hskip\labelsep{\bf Example}]
 \refstepcounter{imn3}{\bf\thesection.\arabic{imn3}.}}%
 {\endtrivlist}

\def\o{\overline}
\makeatletter
\long\def\UR#1{\leavevmode\setbox\@tempboxa\hbox{#1}\@tempdima\fboxrule
    \advance\@tempdima \fboxsep \advance\@tempdima \dp\@tempboxa
   \hbox{\lower \@tempdima\hbox
  {\vbox{\hrule \@height \fboxrule
          \hbox{  \hskip\fboxsep
          \vbox{\vskip\fboxsep \box\@tempboxa\vskip\fboxsep}\hskip
                 \fboxsep\vrule \@width \fboxrule}%
                  }}}}
\makeatother

\makeatletter
\long\def\LR#1{\leavevmode\setbox\@tempboxa\hbox{#1}\@tempdima\fboxrule
    \advance\@tempdima \fboxsep \advance\@tempdima \dp\@tempboxa
   \hbox{\lower \@tempdima\hbox
  {\vbox{
          \hbox{  \hskip\fboxsep
          \vbox{\vskip\fboxsep \box\@tempboxa\vskip\fboxsep}\hskip
                 \fboxsep\vrule \@width \fboxrule}%
                 \hrule \@height \fboxrule}}}}
\makeatother

\makeatletter
\long\def\UL#1{\leavevmode\setbox\@tempboxa\hbox{#1}\@tempdima\fboxrule
    \advance\@tempdima \fboxsep \advance\@tempdima \dp\@tempboxa
   \hbox{\lower \@tempdima\hbox
  {\vbox{\hrule \@height \fboxrule
          \hbox{\vrule \@width \fboxrule \hskip\fboxsep
          \vbox{\vskip\fboxsep \box\@tempboxa\vskip\fboxsep}\hskip
                 \fboxsep }%
                  }}}}
\makeatother

\makeatletter
\long\def\LL#1{\leavevmode\setbox\@tempboxa\hbox{#1}\@tempdima\fboxrule
    \advance\@tempdima \fboxsep \advance\@tempdima \dp\@tempboxa
   \hbox{\lower \@tempdima\hbox
  {\vbox{
          \hbox{\vrule \@width \fboxrule \hskip\fboxsep
          \vbox{\vskip\fboxsep \box\@tempboxa\vskip\fboxsep}\hskip
                 \fboxsep }%
                 \hrule \@height \fboxrule}}}}
\makeatother

\date{}

\title{Unsolved Problems in Virtual Knot Theory and Combinatorial Knot Theory}

\author{
Roger Fenn\\
{\em\footnotesize Department of Mathematics,}\\
{\em\footnotesize University of Sussex, England}\\
{\em\footnotesize rogerf@sussex.ac.uk}\\\\
Denis P.~Ilyutko\\
{\em\footnotesize Department of Mechanics and Mathematics,}\\
{\em\footnotesize Moscow State University, Russia}\\
{\em\footnotesize ilyutko@yandex.ru}\\\\
Louis H.~Kauffman\\
{\em\footnotesize Department of Mathematics, Statistics and
Computer Science,}\\
{\em\footnotesize University of Illinois at Chicago, USA}\\
{\em\footnotesize kauffman@uic.edu}\\\\
Vassily O.~Manturov\\
{\em\footnotesize Department of Fundamental Sciences,}\\
{\em\footnotesize Bauman Moscow State Technical
University, Russia and}\\
{\em\footnotesize Laboratory of Quantum Topology,}\\
{\em\footnotesize Chelyabinsk State University,  Russia}\\
{\em\footnotesize vomanturov@yandex.ru}}

\begin{document}

\renewcommand{\baselinestretch}{0.75}
 \maketitle
\renewcommand{\baselinestretch}{1}

\abstract{This paper is a concise introduction to virtual knot theory, coupled with a list of research problems in this field.}

\newpage
\tableofcontents

\section{Introduction}

The purpose of this paper is to give an introduction to virtual knot
theory and to record a collection of research problems that the
authors have found fascinating. The second section of the paper
introduces the theory and discusses some problems in that context.
Starting from the third sections, we present specific problems. This
paper is an expanded and revised version of our earlier paper on the
subject of problems in virtual knot theory (see~\cite{FKM}). Here we
include a wider selection of problems, including some problems of a
combinatorial flavor in classical knot theory and other problems in
extensions and variants of the theory of knots and links.

We would like to take this opportunity to acknowledge the many
people who have worked on the theory of virtual knots and links.

Researchers in this area explicitly mentioned or referenced in this
paper are: R.\,S.~Avdeev, V.\,G.~Bardakov, A.~Bartholomew,
D.~Bar-Natan, S.~Budden, S.~Carter, K.~Chu, M.\,W.~Chrisman, H.~Dye,
R.~Fenn, R.~Furmaniak, M.~Gousssarov, J.~Green, A.~Henrich,
D.~Hrencecin, D.\,P.~Ilyutko, D.~Jelsovsky, M.~Jordan, T.~Kadokami,
A.~Kaestner, N.~Kamada, S.~Kamada, L.~Kauffman, T.~Kishino,
G.~Kuperberg, S.~Lambropoulou, V.\,O.~Manturov, S.~Nelson,
M.~Polyak, D.\,E.~Radford, S.~Satoh, J.~Sawollek, M.~Saito,
W.~Schellhorn, D.~Silver, V.~Turaev, V.\,V.~Vershinin, O.~Viro,
S.~Williams, P.~Zinn-Justin and J.\,B.~Zuber. See
\cite{AsPrSi,Avdeev,Bardakov,FB,Dror,Blo,FBu,FBu2,Carter,CSH,CSW1,CS,ZCheng,MC_lattice,MC,MC_virtmove,MC_prime,CD_three,
cm_fiber,cm_parity,Chu:FlatKnots,CT,DH,Dye,D1,D2,KDK,DKMin,DK,LinkHtpy,Arrow,DKM,EKT,FENN,FJK,FKM,FRR,FR,FRS2,FRS1,
FRS4,FRS3,FTay,FTur,GPV,HR,HRK,Ily2,IM1,IM2,IM3,IM4,IMN1,IMN2,IS1,KADOKAMI,NKamada,NKamada1,NSKamada,Kamada,KANENOBU,VKT,LKIntro,SVKT,DVK,SLK,NCW1,KD,
ExtBr,NCW2,Affine,KL,VirtualM,Hard,CVBraid,KLTOP,TL,KM1,KaMa,GEN,KS,KiSa,KIS,knotilus,KM,KUP,Mant4,Mant5,Mant6,Mant7,Mant8,Mant9,
Mant10,Mant11,Mant12,Mant13,Mant14,Mant15,Mant16,Mant17,Mant18,Mant19,Mant20,Mant21,Mant22,Mant23,Mant24,Mant25,Mant26,
Mant27,Mant28,Mant29,manturov_compact_long,Mant30,Mant31,
Mant32,Mant33,Mant34,Mant35,Mant37,Mant36,Mant38,Ma23,ManIly,NELSON,
Nelson1,SATOH,SAW,SAW2,Schellhorn,SW,SW1,SW2,SW3,MThist,T,TURAEV,Ver,JZ,JZ1}.
We apologize to anyone who was left out of this list of participant
researchers, and we hope that the  problems described herein will
stimulate people on and off this list to enjoy the beauty of virtual
knot theory!

 \section*{Acknowledgments}
The second and fourth authors were partially supported by grants of
the Russian Government 11.G34.31.0053, RF President NSh --
1410.2012.1, RFBR 12-01-31432, 13-01-00664-a, 13-01-00830-a,
14-01-91161 and 14-01-31288. The fourth author is also partially
supported by Laboratory of Quantum Topology of Che\-lyabinsk State
University (Russian Federation government grant 14.Z50.31.0020). It
gives the third author great pleasure to acknowledge support from
NSF Grant DMS--0245588.


\section{Knot Theories}


\subsection{Virtual knot theory}

\begin{figure}
  \centering\includegraphics[width=12cm]{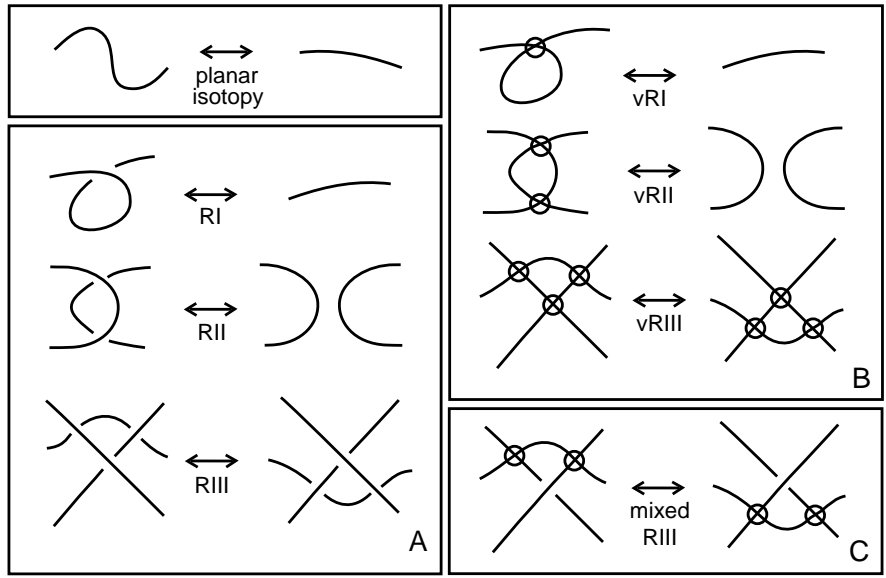}
  \caption{Generalized Reidemeister moves for virtual knots}
  \label{fg:fikm1}
 \end{figure}

Knot theory studies the embeddings of curves in three-dimensional
space. Virtual knot theory studies the  embeddings of curves in
thickened surfaces of arbitrary genus, up to the addition and
removal of empty handles from the surface.  Virtual knots have a
special  diagrammatic theory that makes handling them very similar
to the handling of classical knot diagrams. In fact, this
diagrammatic theory simply involves adding a new type of crossing to
the knot diagrams, a {\em virtual crossing} that is neither under
nor over. From a combinatorial point of view, the virtual crossings
are artifacts of the representation of the virtual knot or link in
the plane. The extension of the Reidemeister moves that takes care
of them respects this viewpoint. A virtual crossing  (see
Fig.~\ref{fg:fikm1}) is represented by two crossing arcs with a
small circle placed around the crossing point.

\begin{figure}
  \centering\includegraphics[width=12cm]{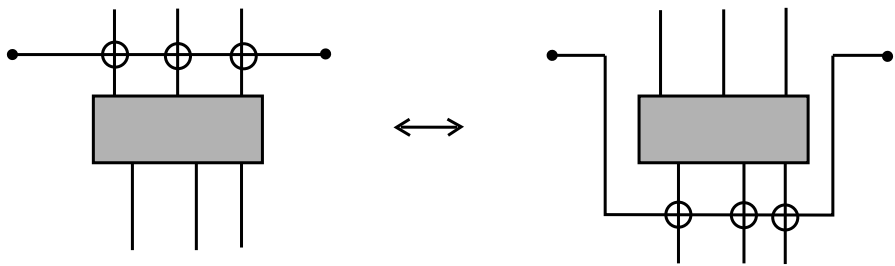}
  \caption{Detour Move}
  \label{fg:fikm2}
 \end{figure}

Moves on virtual diagrams generalize the Reidemeister moves for
classical knot and link diagrams, see Fig.~\ref{fg:fikm1}.  One can
summarize the moves on virtual diagrams by saying that the classical
crossings interact with one another according to the usual
Reidemeister moves. One adds the detour moves for consecutive
sequences of virtual crossings and this completes the description of
the moves on virtual diagrams. It is a consequence of  moves (B) and
(C) in Fig.~\ref{fg:fikm1} that an arc going through any consecutive
sequence of virtual crossings can be moved anywhere in the diagram
keeping the endpoints fixed; the places  where the moved arc crosses
the diagram become new virtual crossings. This replacement is the
{\em detour move}, see Fig.~\ref{fg:fikm2}.

One can generalize many structures in classical knot theory to the
virtual domain, and use the virtual knots to test the limits of
classical problems such as the question whether the Jones
polynomial~\cite{NMR,JO,JO1,JO2,NKamada,NKamada1,K,KA89,QJones,Mur,WITT}
detects knots and other classical problems. Counterexamples often
exist in the virtual domain, and it is an open problem whether these
counterexamples are equivalent (by addition and subtraction of empty
handles) to classical knots and links. Virtual knot theory is a
significant domain to be investigated for its own sake and for a
deeper understanding of classical knot theory.

Another way to understand the meaning of virtual diagrams is to
regard them as representatives for oriented Gauss codes (Gauss
diagrams)~\cite{GPV,GAUSS,VKT}. Such codes do not always have planar
realizations and an attempt to embed such a code in the plane leads
to the production of the virtual crossings. The detour move makes
the particular choice of virtual crossings irrelevant. Virtual
equivalence is the same as the equivalence relation generated on the
collection of oriented Gauss codes modulo an abstract set of
Reidemeister moves on the codes.

 \subsection{Flat virtual knots and links}

Every classical knot or link diagram can be regarded as a
$4$-regular plane graph with extra structure at the nodes. This
extra structure is usually indicated by the over and under crossing
conventions that give instructions for constructing an embedding of
the link in three dimensional space from the diagram.  If we take
the diagram without this extra structure, it is the shadow of some
link in three dimensional space, but the weaving of that link is not
specified. It is well known that if one is allowed to apply the
Reidemeister moves to such a shadow (without regard to the types of
crossing since they are not specified) then the shadow can be
reduced to a disjoint union of circles. This reduction is no longer
true for virtual links. More precisely, let a {\em flat virtual
diagram} be a diagram with virtual crossings as we have described
them and {\em flat crossings} consisting in undecorated nodes of the
$4$-regular plane graph. Virtual crossings are flat crossings that
have been decorated by a small circle. Two flat virtual diagrams are
{\em equivalent} if there is a sequence of generalized flat
Reidemeister moves (as illustrated in Fig.~\ref{fg:fikm1}) taking
one to the other. A generalized flat Reidemeister move is any move
as shown in Fig.~\ref{fg:fikm1}, but one can ignore the over or
under crossing structure. Note that in studying flat virtual knots
the rules for changing virtual crossings among themselves and the
rules for changing flat crossings among themselves are identical.
However, detour moves as in Fig.~\ref{fg:fikm1}C are available for
virtual crossings with respect to flat crossings and not the other
way around.

We shall say that a virtual diagram {\em overlies} a flat diagram if
the virtual diagram is obtained from the flat diagram by choosing a
crossing type for each flat crossing in the virtual diagram. To each
virtual diagram $K$ there is an associated flat diagram $F(K)$ that
is obtained by forgetting the extra structure at the classical
crossings in $K.$ Note that if $K$ is equivalent to $K'$ as virtual
diagrams, then $F(K)$ is equivalent to $F(K')$ as flat virtual
diagrams. Thus, if we can show that $F(K)$ is not reducible to a
disjoint union of circles, then it will follow that $K$ is a
non-trivial virtual link.

 \begin{figure}
  \centering\includegraphics[width=10cm]{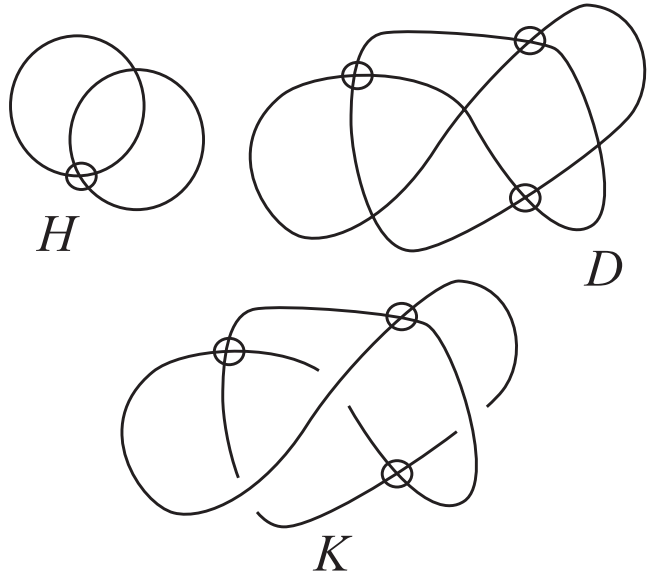}
  \caption{Flats $H$ and $D$, and the knot $K$}\label{fg:fikm3}
 \end{figure}

Fig.~\ref{fg:fikm3} illustrates an example of a flat virtual link
$H.$ This link cannot be undone in the flat category because it has
an odd number of virtual crossings between its two components and
each generalized Reidemeister move preserves the parity of the
number of virtual crossings between components.  Also illustrated in
Fig.~\ref{fg:fikm3} is a flat diagram $D$ and a virtual knot $K$
that overlies it. This example is given in~\cite{VKT}. The knot
shown is undetectable by many invariants (fundamental group, Jones
polynomial) but it is knotted.  The flat virtual diagrams present a
challenge for the construction of new invariants.  It is important
to understand the structure of flat virtual knots and links. This
structure lies at the heart of the comparison of classical and
virtual links.  Simpler and more fundamental than flat virtual knots
and links are the {\em free knots and links}~\cite{Mant31,Mant32}
corresponding to Gauss diagrams with no signs and no arrows. We will
list a number of problems about free knots below.

 \subsection{Interpretation of virtual knots as stable classes of links in  thickened surfaces}

There is a useful topological interpretation for this virtual theory
in terms of embeddings of links in thickened surfaces,
see~\cite{VKT,DVK,KUP}. Regard each virtual crossing as a shorthand
for a detour of one of the arcs in the crossing through a $1$-handle
that has been attached to the 2-sphere of the original diagram. By
interpreting each virtual crossing in this way, we obtain an
embedding of a collection of circles into a thickened surface $S_{g}
\times \mathbb{R}$, where $g$ is the number of virtual crossings in
the original diagram $L$, $S_{g}$ is a compact oriented surface of
genus $g$ and $\mathbb{R}$ denotes the real line.  We say that two
such surface embeddings are {\em stably equivalent} if one can be
obtained from another by isotopy in the thickened surfaces,
homeomorphisms of the surfaces and the addition or subtraction of
empty handles.  Then we have the following theorem.

  \begin{theorem}[see~\cite{VKT,DKT,KUP}]
Two virtual link diagrams are equivalent if and only if their
correspondent surface embeddings are stably equivalent.
 \end{theorem}

Here long knots (or, equivalently $1-1$ tangles) come into play.
Having a knot, we can break it at some point and take its ends to
infinity (say, in a way that they coincide with the horizontal axis
line in the plane). One can study isotopy classes of such knots. A
well-known theorem says that in the classical case, knot theory
coincides with long knot theory. However, this is not the case for
virtual knots. By breaking the same virtual knot at different
points, one can obtain non-isotopic long knots~\cite{BFKK,Mant20}.
Furthermore, even if the initial knot is trivial, the resulting long
knot may not be trivial. A  connected sum of two trivial virtual
diagrams may not be trivial in the compact case. The phenomenon
occurs because these two knot diagrams may  be non-trivial in the
long category. It is sometimes more convenient to consider long
virtual knots rather than compact virtual knots, since connected sum
is well-defined for long knots.

Unlike classical knots, the connected sum of long virtual knots is
not commutative~\cite{BFKK,Mant19,Mant20}. Thus, if we show that two
long knots $K_{1}$ and $K_{2}$ do not commute, then we see that they
are different and both non-classical.

\begin{figure}
  \begin{center}
  {\tt    \setlength{\unitlength}{0.92pt}
  \begin{picture}(195,239)
  \thicklines   \put(111,11){\line(0,1){47}}
              \put(114,99){\line(0,1){47}}
              \put(95,10){\line(0,1){47}}
              \put(92,99){\line(0,1){47}}
              \put(111,58){\line(1,0){29}}
              \put(114,98){\line(1,0){29}}
              \put(65,58){\line(1,0){29}}
              \put(62,99){\line(1,0){29}}
              \put(24,59){\line(0,-1){12}}
              \put(163,79){\circle{20}}
              \put(45,78){\circle{20}}
              \put(24,60){\line(1,1){39}}
              \put(183,45){\line(0,1){17}}
              \put(133,43){\line(1,0){50}}
              \put(133,51){\line(0,-1){8}}
              \put(133,95){\line(0,-1){34}}
              \put(133,114){\line(0,-1){10}}
              \put(184,114){\line(-1,0){49}}
              \put(183,100){\line(0,1){13}}
              \put(77,46){\line(0,1){7}}
              \put(24,45){\line(1,0){53}}
              \put(77,95){\line(0,-1){33}}
              \put(77,113){\line(0,-1){10}}
              \put(25,113){\line(1,0){52}}
              \put(25,99){\line(0,1){13}}
              \put(141,59){\line(1,1){40}}
              \put(145,98){\line(1,-1){37}}
              \put(24,99){\line(1,-1){40}}
              \put(12,212){\line(1,-1){40}}
              \put(53,213){\line(1,0){41}}
              \put(96,211){\line(1,-1){37}}
              \put(52,171){\line(1,0){40}}
              \put(92,172){\line(1,1){40}}
              \put(13,212){\line(0,1){13}}
              \put(13,226){\line(1,0){52}}
              \put(65,226){\line(0,-1){10}}
              \put(65,208){\line(0,-1){33}}
              \put(12,171){\line(0,-1){12}}
              \put(12,158){\line(1,0){53}}
              \put(65,159){\line(0,1){7}}
              \put(134,213){\line(0,1){13}}
              \put(135,227){\line(-1,0){49}}
              \put(84,227){\line(0,-1){10}}
              \put(84,208){\line(0,-1){34}}
              \put(84,164){\line(0,-1){8}}
              \put(84,156){\line(1,0){50}}
              \put(134,158){\line(0,1){17}}
              \put(12,173){\line(1,1){39}}
              \put(33,191){\circle{20}}
              \put(114,192){\circle{20}}
              \put(139,175){\makebox(32,31){$K$}}
  \end{picture}}
  \end{center}
 \caption{Kishino and parts}\label{fg:fikm4}
 \end{figure}
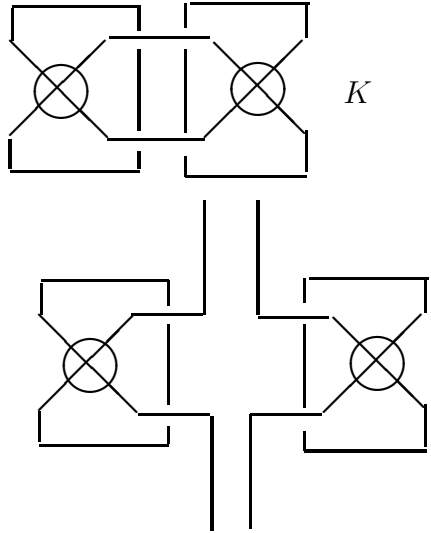

A typical example of such knots is the two parts of the Kishino
knot, see Fig.~\ref{fg:fikm4}.

We have a natural map
 $$
\langle\mbox{Long virtual knots }\rangle\to \langle\mbox{Oriented
compact virtual knots}\rangle,
 $$
obtained by taking two infinite ends of the long knots together to
make a compact knot. This map is obviously well defined.

Note that when the parts of the Kishino knot are closed they become
unknots.

This map allows one to construct long virtual knot invariants from
classical invariants, i.e., just to regard compact knot invariants
as long knot invariants. There is no well-defined inverse for this
map. The long category can also be applied for the case of flat
virtual knots, where all problems formulated above occur as well.


\section{Switching and Virtualizing}

\begin{figure}
  \begin{center}
{\tt    \setlength{\unitlength}{0.92pt}
\begin{picture}(326,162)
\thinlines    \put(190,141){\line(0,-1){82}}
              \put(242,34){\makebox(41,41){$s(i)$}}
              \put(170,10){\makebox(41,41){$v(i)$}}
              \put(1,35){\makebox(41,42){$i$}}
              \put(292,81){\line(1,0){33}}
              \put(246,81){\line(1,0){35}}
              \put(286,161){\line(0,-1){160}}
              \put(143,81){\circle{16}}
              \put(191,80){\circle{16}}
              \put(190,59){\line(-1,0){25}}
              \put(167,22){\line(0,-1){20}}
              \put(144,22){\line(1,0){22}}
              \put(143,102){\line(0,-1){80}}
              \put(165,103){\line(-1,0){22}}
              \put(165,87){\line(0,1){15}}
              \put(165,59){\line(0,1){17}}
              \put(167,141){\line(1,0){23}}
              \put(166,161){\line(0,-1){20}}
              \put(126,81){\line(1,0){80}}
              \put(46,74){\line(0,-1){73}}
              \put(46,161){\line(0,-1){74}}
              \put(6,81){\line(1,0){80}}
  \end{picture}}
  \end{center}
  \caption{Switching and virtualizing a crossing}\label{fg:fikm5}
  \end{figure}
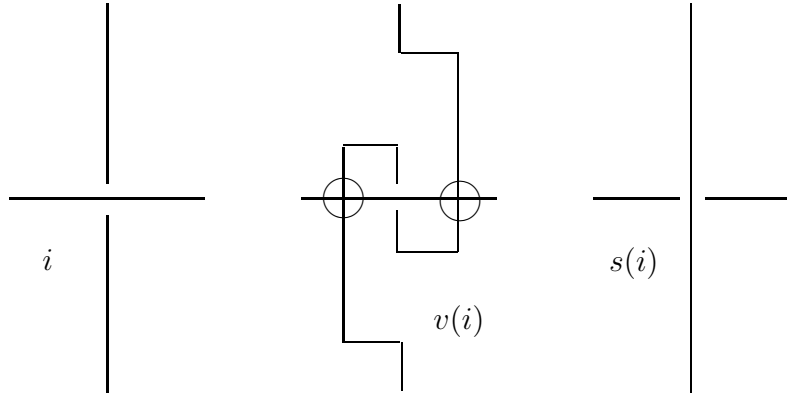

Given a crossing $i$ in a link diagram, we define $s(i)$ to be the
result of {\em switching} that crossing so that the undercrossing
arc becomes an overcrossing arc and vice versa. We also define the
{\em virtualization} $v(i)$ of the crossing by the local replacement
indicated in Fig.~\ref{fg:fikm5}. In this figure we illustrate how
in the virtualization of the crossing the  original crossing is
replaced by a crossing that is flanked by two virtual crossings.

Suppose that $K$ is a (virtual or classical) diagram with a
classical crossing labeled $i.$  Let $K^{v(i)}$ be the diagram
obtained from $K$ by virtualizing the crossing $i$ while leaving the
rest of the diagram just as before. Let $K^{s(i)}$ be the diagram
obtained from $K$ by switching the crossing $i$ while leaving the
rest of the diagram just as before. Then it follows directly from
the definition of the Jones polynomial that $$V_{K^{s(i)}}(t) =
V_{K^{v(i)}}(t).$$ \noindent As far as the Jones polynomial is
concerned, switching a crossing and virtualizing a crossing look the
same.

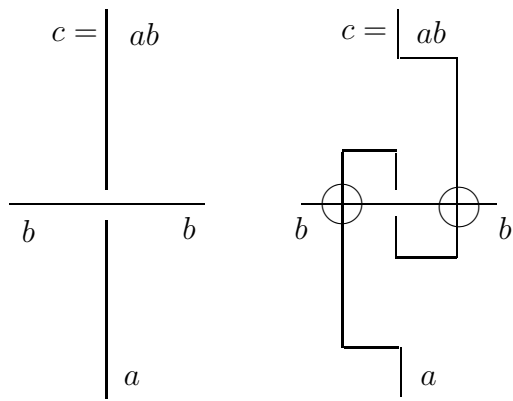
\begin{figure}
  \begin{center}
{\tt    \setlength{\unitlength}{0.92pt}
\begin{picture}(214,164)
\thinlines    \put(197,64){\makebox(16,16){$b$}}
              \put(136,144){\makebox(22,18){$c=$}}
              \put(17,143){\makebox(22,18){$c=$}}
              \put(165,143){\makebox(19,20){$ab$}}
              \put(113,64){\makebox(16,16){$b$}}
              \put(67,64){\makebox(16,16){$b$}}
              \put(43,1){\makebox(17,19){$a$}}
              \put(47,142){\makebox(19,20){$ab$}}
              \put(1,63){\makebox(16,16){$b$}}
              \put(165,1){\makebox(17,19){$a$}}
              \put(185,142){\line(0,-1){82}}
              \put(138,82){\circle{16}}
              \put(186,81){\circle{16}}
              \put(185,60){\line(-1,0){25}}
              \put(162,23){\line(0,-1){20}}
              \put(139,23){\line(1,0){22}}
              \put(138,103){\line(0,-1){80}}
              \put(160,104){\line(-1,0){22}}
              \put(160,88){\line(0,1){15}}
              \put(160,60){\line(0,1){17}}
              \put(162,142){\line(1,0){23}}
              \put(161,162){\line(0,-1){20}}
              \put(121,82){\line(1,0){80}}
              \put(41,75){\line(0,-1){73}}
              \put(41,162){\line(0,-1){74}}
              \put(1,82){\line(1,0){80}}
\end{picture}}
 \end{center}
 \caption{$IQ(Virt(K)) = IQ(K)$}\label{fg:fikm6}
 \end{figure}

The involutory quandle~\cite{KNOTS} is an algebraic invariant
equivalent to the fundamental group of the double branched cover of
a knot or link in the classical case. In this algebraic system one
associates a generator of the algebra $IQ(K)$ to each arc of the
diagram $K$ and there is a relation of the form $c = ab$ at each
crossing, where $ab$ denotes the (non-associative) algebra product
of $a$ and $b$ in $IQ(K).$ See Fig.~\ref{fg:fikm6}. In this Figure
we have illustrated through the local relations the fact that
 $$
IQ(K^{v(i)}) = IQ(K).
 $$
As far the involutory quandle is concerned, the original crossing
and the virtualized crossing look the same.

If a classical knot is actually knotted, then its involutory quandle
is non-trivial~\cite{W}. Hence if we start with a non-trivial
classical knot, we can virtualize any subset of its crossings to
obtain a virtual knot that is still non-trivial. There is a subset
$A$ of the crossings of a classical knot $K$ such that the knot $SK$
obtained by switching these crossings is an unknot.  Let $Virt(K)$
denote the virtual diagram obtained from $A$ by virtualizing the
crossings in the subset $A.$  By the above discussion the Jones
polynomial of $Virt(K)$ is the same as the Jones polynomial of $SK$,
and this is $1$ since $SK$ is unknotted. On the other hand, the $IQ$
of $Virt(K)$ is the same as the $IQ$ of $K$, and hence if $K$ is
knotted, then so is $Virt(K).$   We have shown that $Virt(K)$ is a
non-trivial virtual knot with unit Jones polynomial.  This completes
the proof of the theorem.

If there exists a classical knot with unit Jones polynomial, then
one of the knots $Virt(K)$ produced by this theorem may be
equivalent to a classical knot.  It is an intricate task to verify
that specific examples of $Virt(K)$ are not classical.

A very fruitful line of new invariants comes about by examining a
generalization of the fundamental group or quandle that we call the
{\em biquandle} of the virtual knot. The biquandle is discussed in
the next Section. Invariants of flat knots (when one has them) are
useful in this regard. If we can verify that the flat knot
$F(Virt(K))$ is non-trivial, then $Virt(K)$ is non-classical. In
this way the search for classical knots with unit Jones polynomial
(see~\cite{MThist} for links) expands to the exploration of the
structure of the infinite collection of virtual knots with unit
Jones polynomial.

Another way of putting this theorem is as follows: In the arena of
knots in thickened surfaces there are many examples of knots with
unit Jones polynomial. Might one of these be equivalent via handle
stabilization to a classical knot? In~\cite{KUP} Kuperberg shows the
uniqueness of the  embedding of minimal genus in the stable class
for a given virtual link. The minimal embedding genus can be
strictly less than the number of virtual crossings in a diagram for
the link.  There are many problems associated with this phenomenon.


\section{Atoms}\label{sec:atoms}


An {\em atom} is a pair: $(M^{2},\Gamma)$ where $M^{2}$ is a closed
$2$-manifold and $\Gamma$ is a 4-valent graph in $M^{2}$ dividing
$M^{2}$ into cells such that these cells admit a checkerboard
coloring (the coloring is also fixed). $\Gamma$ is called the {\em
frame} of the atom,
see~\cite{Fom2,Fom3,Fom4,Fom,Mant1,Mant3,Mant19,Mant21}.

Atoms are considered up to natural equivalence, that is, up to
homeomorphisms of the underlying manifold $M^{2}$ mapping the frame
to the frame and black cells to black cells. From this point of
view, an atom can be recovered from the frame together with the
following combinatorial structure:
 \begin{enumerate}
  \item
{\em $A$-structure}: This indicates which edges for each vertex are
{\em opposite edges}. That is, it indicates the cyclic structure at
the vertex.
  \item
{\em $B$-structure}: This indicates pairs of ``black angles''. That
is, one divides the four edges emanating from each vertex into two
sets of adjacent (not opposite) edges such that the black cells are
locally attached along these pairs of adjacent edges.
 \end{enumerate}

Given a virtual knot diagram which has its regions black and white
colored like a chess board, one can construct the corresponding atom
as follows. Classical crossings correspond to the vertices of the
atom, and generate both the $A$-structure and the $B$-structure at
these vertices (the $B$-structure comes from over/under
information). Thus, an atom is uniquely determined by a virtual knot
diagram. It is easy to see that the inverse operation is well
defined modulo virtualization. Thus the atom knows everything about
the bracket polynomial (Jones polynomial) of the virtual link.

The crucial notions here are the minimal genus of the atom and the
orientability of the atom. For instance, for each link diagram with
a corresponding orientable atom (all classical link diagrams are in
this class), all degrees of the bracket are congruent modulo four
while in the non-orientable case they are congruent only modulo two.

The orientability condition is crucial in the construction of the
Khovanov homology theory for virtual links as
in~\cite{Mant19,Mant21}.


 \section{Biquandles}


 \subsection{Main constructions}

In this section we give a sketch of some recent approaches to
invariants of virtual knots and links.

\begin{figure}
  \begin{center}
{\tt    \setlength{\unitlength}{0.92pt}
\begin{picture}(280,163)
\thicklines   \put(131,162){\line(0,-1){159}}
              \put(1,1){\framebox(278,161){}}
              \put(231,104){\makebox(21,19){$a^{\o{b}} = a\, \UL{b}$}}
              \put(167,59){\makebox(20,15){$b_{\o{a}} = b\, \LL{a}$}}
              \put(233,62){\makebox(20,16){$b$}}
              \put(212,38){\makebox(18,17){$a$}}
              \put(73,106){\makebox(23,23){$a^{b} = a\, \UR{b}$}}
              \put(77,56){\makebox(32,22){$b_{a} = b\, \LR{a}$}}
              \put(10,86){\makebox(20,16){$b$}}
              \put(52,43){\makebox(18,17){$a$}}
              \put(211,89){\vector(0,1){33}}
              \put(211,43){\vector(0,1){32}}
              \put(250,82){\vector(-1,0){78}}
              \put(51,90){\vector(0,1){33}}
              \put(51,43){\vector(0,1){34}}
              \put(11,83){\vector(1,0){80}}
 \end{picture}}
 \end{center}
 \caption{Biquandle relations at a crossing}\label{fg:fikm7}
 \end{figure}
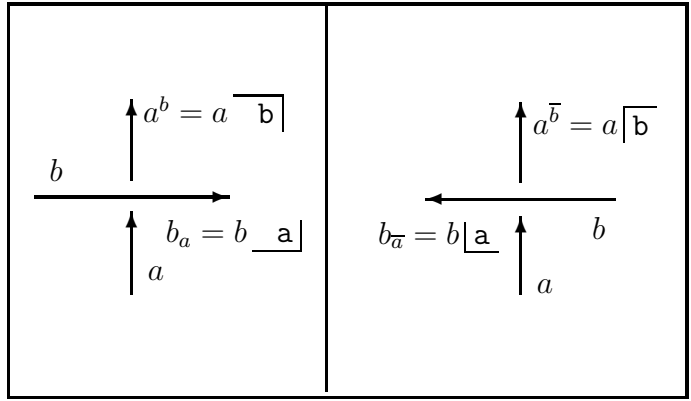

A {\em biquandle}~\cite{FB,FBu,CS,FJK,DVK,GEN} is an algebra with
four binary operations written
$a^b,\,a_b,\,a^{\overline{b}},\,a_{\overline{b}}$ together with some
relations which we will indicate below. The {\em fundamental}
biquandle is associated with a link diagram and is invariant under
the generalized Reidemeister moves for virtual knots and links.  The
operations in this algebra are motivated by the formation of labels
for the edges of the diagram, see Fig~\ref{fg:fikm7}. In this figure
we have shown the format for the operations in a biquandle. The
overcrossing arc has two labels, one on each side of the crossing.
There is an algebra element labeling each {\em edge} of the diagram.
An edge of the diagram corresponds to an edge of the underlying
plane graph of that diagram.

Let the edges oriented toward a crossing in a diagram be called the
{\em input} edges for the crossing, and the edges oriented away from
the crossing be called the {\em output} edges for the crossing. Let
$a$ and $b$ be the input edges for a positive crossing, with $a$ the
label of the undercrossing input and $b$ the label on the
overcrossing input. In the biquandle, we label the undercrossing
output by $$c=a^{b},$$ while the overcrossing output is labeled
 $$d= b_{a}.$$

The labelling for the negative crossing is similar using the other
two operations.

To form the fundamental biquandle, $BQ(K)$, we take one generator
for each edge of the diagram and two relations at each crossing (as
described above).

Another way to write this formalism for the biquandle is as follows
 $$
a^{b} = a\, \UR{b},\, a_{b} = a\, \LR{b},\, a^{\o{b}} = a\,
\UL{b},\, a_{\o{b}} = a\, \LL{b}.
 $$
We call this the {\em operator formalism} for the biquandle.

These considerations lead to the following definition.

 \begin{definition}
A {\em biquandle} $B$ is a set with four binary operations indicated
above:  $a^{b} \,\mbox{,} \, a^{\o{b}} \, \mbox{,} \,  a_{b}
\,\mbox{,} \, a_{\o{b}}.$ We shall refer to the operations with
barred variables as the {\em left} operations and the operations
without barred variables as the {\em right} operations. The
biquandle is closed under these operations and the following axioms
are satisfied:
 \begin{enumerate}
  \item
Given an element $a$ in $B$, then there exists an $x$ in the
biquandle such that $x=a_{x}$ and $a = x^{a}.$ There also exists a
$y$ in the biquandle such that $y=a^{\o{y}}$ and $a = y_{\o{a}}.$
  \item
For any elements $a$ and $b$ in $B$ we have
 $$
a = a^{b \o{b_{a}}},\,  b= b_{a \o{a^{b}}},\,a =
a^{\o{b}b_{\o{a}}},\, b= b_{\o{a} a^{\o{b}}}.
 $$
 \item
Given elements $a$ and $b$ in $B$ then there exist elements
$x,\,y,\,z,\,t$ such that $x_{b}=a$, $y^{\overline{a}}= b$, $b^x=y$,
$a_{\overline{y}}=x$ and $t^a=b$, $a_t=z$, $z_{\overline{b}}=a$,
$b^{\overline{z}}=t$.

The biquandle is called {\em strong} if $x,\,y,\,z,\,t$ are uniquely
defined and we then write $x=a_{b^{-1}}, y=b^{\overline{a}^{-1}},
t=b^{a^{-1}}, z=a_{\overline{b}^{-1}}$, reflecting the invertive
nature of the elements.
  \item
For any $a$, $b$, $c$ in $B$ the following equations hold and the
same equations hold when all right operations are replaced in these
equations by left operations:
 $$
a^{b c} = a^{c_{b} b^{c}},\, c_{b a} = c_{a^{b} b_{a}},\,
(b_{a})^{c_{a^{b}}} = (b^{c})_{a^{c_{b}}}.
 $$
 \end{enumerate}
 \end{definition}

These axioms are transcriptions of the Reidemeister moves.The first
axiom transcribes the first Reidemeister move. The second axiom
transcribes the directly oriented second Reidemeister move. The
third axiom transcribes the reverse oriented Reidemeister move. The
fourth axiom transcribes the third Reidemeister move. Much more work
is needed in exploring these  algebras and their applications to
knot theory.

We may simplify the appearance of these conditions by defining
 $$
S(a,b)=(b_a,a^b),\quad
\overline{S}(a,b)=(b^{\overline{a}},a_{\overline{b}})$$ and in the
case of a strong
biquandle,$$S^+_-(a,b)=(b^{a_{b^{-1}}},a_{b^{-1}}),\quad
S^-_+(a,b)=(b^{a^{-1}},a_{b^{a^{-1}}})$$ and
$${\overline{S}}^{\lower.5ex\hbox{\ $\scriptstyle+$}}_{\ -}(a,b)=
(b_{\ \overline{a^{\overline{b}^{-1}}}}\ , \ a^{\overline{b}^{-1}})=
(b_{\ \overline{a^{b_{a^{-1}}}}} \ , \ a^{b_{a^{-1}}})$$ \ and \
$${\overline{S}}^{\lower.5ex\hbox{\ $\scriptstyle-$}}_{\ +}(a,b)=
(b_{\ \overline{a}^{-1}} \ , \
a^{\overline{b_{\overline{a}^{-1}}}})= (b_{{a}^{b^{-1}}} \ , \
a^{\overline{b_{{a}^{b^{-1}}}}})
 $$
which we call the {\em sideways} operators. The conditions then
reduce to $$S\overline{S}=\overline{S}S=1,$$ $$ (S\times 1) (1\times
S) (S\times 1) = (1\times S) (S\times 1) (1\times S)$$
$$\overline{S}^-_+S^+_-=S^-_+\overline{S}^+_-=1$$ and finally all
the sideways operators leave the diagonal
 $$
\Delta=\{(a,a)|a\in X\}
 $$
invariant. There is a different and possibly simpler approach
in~\cite{Ftref}.

Here the sideways operator is used to define the up and down
actions. So if $F\colon X^2\to X^2$ denotes the sideways operator
then let
 $$
F(a,b)=(b_a,a^b)=\bigl(f_a(b),f^b(a)\bigr).
 $$
The sideways operator is required to satisfy the following:
 \begin{enumerate}
  \item[1)]
$F$ is a bijection,
  \item[2)]
$f_*$ is a bijection,
  \item[3)]
$f^*$ is a bijection.
  \end{enumerate}
Then $S$, the switch operator, is defined by
  $$
S(b_a,a)=(a^b,b).
 $$
It turns out, symmetrically, that
 \begin{enumerate}
  \item[1)]
$S$ is a bijection,
  \item[2)]
$s_*$ is a bijection,
  \item[3)]
$s^*$ is a bijection,
  \end{enumerate}
where $S(x,y)=\bigl(s^x(y),s_y(x)\bigr).$

The advantage of this approach is that the conditions for $S$ to
satisfy the Yang--Baxter type equations above is that
 $$
1)\ a^{bc_b}=a^{cb^c},\quad 2)\ c_{ba^b}=c_{ab_a}\quad
\hbox{and}\quad 3)\ {b_a}^{c_a}={b^c}_{a^c}
 $$
and for a biquandle $a^a=a_a$, which are considerably simpler than
the above.

\subsection{The Alexander biquandle}

It is not hard to see that the following equations in a module over
$\mathbb{Z}[s,s^{-1},t,t^{-1}]$ give a biquandle structure:
 \begin{gather*}
a^b=a\,\UR{b} = ta + (1-st)b,\quad a_b=a\,\LR{b} = sa\\
a^{\overline{b}}=a\,\UL{b} = t^{-1}a + (1-s^{-1}t^{-1})b,\quad
a_{\overline{b}}=a\,\LL{b} = s^{-1}a.
 \end{gather*}
We shall refer to this structure, with the equations given above, as
the {\em Alexander Biquandle}.

Just as one can define the Alexander Module of a classical knot, we
have the Alexander Biquandle of a virtual knot or link, obtained by
taking one generator for each {\em edge} of the projected graph of
the  knot diagram and taking the module relations in the above
linear form. Let $ABQ(K)$ denote this module structure for an
oriented link $K$. That is, $ABQ(K)$ is the module generated by the
edges of the diagram, factored by the submodule generated by the
relations. This module then has a biquandle structure specified by
the operations defined above for an Alexander Biquandle.

The determinant of the matrix of relations obtained from the
crossings of a diagram gives a polynomial invariant (up to
multiplication by $\pm s^{i}t^{j}$ for integers $i$ and $j$) of
knots and links that we denote by $G_{K}(s,t)$ and call the {\em
generalized Alexander polynomial}. {\it This polynomial vanishes on
classical knots, but is remarkably successful at detecting virtual
knots and links.} In fact $G_{K}(s,t)$ is the same as the polynomial
invariant of virtual knots of Sawollek~\cite{SAW} and defined by an
alternative method by Silver and Williams~\cite{SW} and by yet
another method by Manturov~\cite{Mant11}. It is a reformulation of
the invariant for knots in surfaces due to the principal
investigator, Jaeger and Saleur~\cite{JKS,KS}.

\begin{figure}
\begin{center} {\tt    \setlength{\unitlength}{0.92pt}
\begin{picture}(264,109)
\thicklines   \put(186,1){\makebox(29,29){$KI$}}
              \put(50,1){\makebox(32,31){$K$}}
              \put(241,70){\circle{20}}
              \put(159,70){\circle{20}}
              \put(105,71){\circle{20}}
              \put(24,70){\circle{20}}
              \put(219,49){\line(-1,0){9}}
              \put(223,89){\line(-1,0){13}}
              \put(180,50){\line(1,0){23}}
              \put(179,90){\line(1,0){24}}
              \put(205,35){\line(1,0){55}}
              \put(206,105){\line(0,-1){70}}
              \put(261,105){\line(-1,0){55}}
              \put(139,90){\line(1,-1){40}}
              \put(223,89){\line(1,-1){37}}
              \put(219,50){\line(1,1){40}}
              \put(140,90){\line(0,1){13}}
              \put(140,104){\line(1,0){52}}
              \put(192,104){\line(0,-1){10}}
              \put(192,86){\line(0,-1){33}}
              \put(139,49){\line(0,-1){12}}
              \put(139,36){\line(1,0){53}}
              \put(192,37){\line(0,1){7}}
              \put(261,91){\line(0,1){13}}
              \put(259,36){\line(0,1){17}}
              \put(139,51){\line(1,1){39}}
              \put(3,52){\line(1,1){39}}
              \put(125,37){\line(0,1){17}}
              \put(75,35){\line(1,0){50}}
              \put(75,43){\line(0,-1){8}}
              \put(75,87){\line(0,-1){34}}
              \put(75,106){\line(0,-1){10}}
              \put(126,106){\line(-1,0){49}}
              \put(125,92){\line(0,1){13}}
              \put(56,38){\line(0,1){7}}
              \put(3,37){\line(1,0){53}}
              \put(3,50){\line(0,-1){12}}
              \put(56,87){\line(0,-1){33}}
              \put(56,105){\line(0,-1){10}}
              \put(4,105){\line(1,0){52}}
              \put(4,91){\line(0,1){13}}
              \put(83,51){\line(1,1){40}}
              \put(43,50){\line(1,0){40}}
              \put(87,90){\line(1,-1){37}}
              \put(44,92){\line(1,0){41}}
              \put(3,91){\line(1,-1){40}}
\end{picture}}
\end{center}
\caption{The knot $K$ and the Kishino diagram $KI$ }\label{fg:fikm8}
 \end{figure}
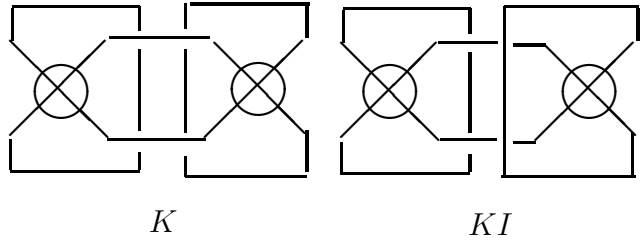

We end this discussion of the Alexander Biquandle with two examples
that show clearly its limitations. View Figure~\ref{fg:fikm8}. In
this figure we illustrate two diagrams labeled $K$ and $KI.$ It is
not hard to calculate that both $G_{K}(s,t)$ and $G_{KI}(s,t)$ are
equal to zero. However, The Alexander Biquandle of $K$ is
non-trivial -- it is isomorphic to the free module over $Z[s,
s^{-1},t, t^{-1}]$ generated by elements $a$ and $b$ subject to the
relation $(s^{-1} - t -1)(a-b) =0.$ Thus $K$ represents a
non-trivial virtual knot. This shows that it is possible for a
non-trivial virtual diagram to be a connected sum of two trivial
virtual diagrams. However, the diagram $KI$ has a trivial Alexander
Biquandle.  In fact the diagram $KI$, discovered by
Kishino~\cite{AsPrSi}, is now known to be knotted and its general
biquandle is non-trivial. The Kishino diagram has been shown
non-trivial by a calculation of the three-strand Jones
polynomial~\cite{KiSa},  by the surface bracket polynomial of Dye
and Kauffman~\cite{Dye,DKMin}, by the $\Xi$-polynomial (the surface
generalization of the Jones polynomial of Manturov~\cite{Mant12},
and its biquandle has been shown to be non-trivial by a quaternionic
biquandle representation~\cite{FB} which we will now briefly
describe.

The quaternionic biquandle is defined by the following operations
where $i^2 = j^2 = k^{2} = ijk = -1$, $ij = -ji = k$, $jk = -kj =
i$, $ki = -ik =j$ in the associative, non-commutative algebra of the
quaternions. The elements $a,\,b$ are in a module over the ring of
integer quaternions:
 \begin{gather*}
a^b=a \UR{b}=j\cdot a+ (1+i)\cdot b,\\
a_b=a \LR{b}=-j\cdot a+ (1+i)\cdot b,\\
a^{\overline{b}}=a \UL{b}=j\cdot a+ (1-i) \cdot b,\\
a_{\overline{b}}=a \LL{b}=-j\cdot a+ (1-i)\cdot b.
 \end{gather*}
Amazingly, one can verify that these operations satisfy the axioms
for the biquandle.

Equivalently, referring back to the previous section, define the
linear biquandle by
 $$
S=\left(\begin{array}{cc}1+i&jt\\ -jt^{-1}&1+i\end{array}\right),
 $$
where $i,\,j$ have their usual meanings as quaternions and $t$ is a
central variable. Let $R$ denote the ring which they determine. Then
as in the Alexander case considered above, for each diagram there is
a square presentation of an $R$-module. We can take the (Study)
determinant of the presentation matrix. In the case of the Kishino
knot this is zero. However the greatest common divisor of the
codimension $1$ determinants is $2+5t^2+2t^4$ showing that this knot
is not classical.

There are many other developments in the theory of virtual knots. We
shall refer to them as need be in composing the problems. We refer
the reader to our list of papers and particularly the introductory
paper~\cite{SVKT} and the book~\cite{ManIly}.


\section{Biracks, Biquandles etc:}


\subsection{Main definitions}

Let $X$ be a set of {\em labels}. Numerous knot invariants can be
defined by a {\em switch} map, $S\colon X^2\to X^2$.

There are various axioms which $S$ may or may not satisfy. They are
listed below in no particular order:
 \begin{enumerate}
  \item
{\em Invertibility of $S$}: The map $S$ is a bijection.
  \item
{\em Right Invertibility of the Binary Products}: Let
$i_*,\,j_*\colon X\to X^2$ be the inclusions $i_a(y)=(a,y)$ and
$j_b(x)=(x,b)$. Let $p,\,q\colon X^2\to X$ be the projections
$p(a,b)=a$ and $q(a,b)=b$. Then the compositions
 $$
X\buildrel i_*\over \longrightarrow X^2\buildrel q\over
\longrightarrow X \hbox{ and } X\buildrel j_*\over \longrightarrow
X^2\buildrel p\over \longrightarrow X
 $$
are bijections.

The first map defines a binary operation, $(a,b)\to a^b$, and the
{\it inverse} of the second defines another, $(a,b)\to b_a$. So the
{switch} map, $S\colon X^2\to X^2$,  is now defined by
 $$
S(a,b_a)=(b,a^b).
 $$

The advantage of the suffix and superfix notation is that brackets
can often be dispensed with: so for example
 $$
a^{bc}=(a^b)^c,\, a^{b_c}=a^{(b_c)},\, {a^b}_c=(a^b)_c,\, \hbox{
etc. }
 $$
On the other hand expressions such as $a^b_c$ are ambiguous and are
not used.

The operations are right invertible. So there are inverse operations
$(a,b)\to a^{b^{-1}}$ and $(a,b)\to a_{b^{-1}}$ such that
 $$
a^{bb^{-1}}= a^{b^{-1}b}=a\,\hbox{ and }\,a_{bb^{-1}}=
a_{b^{-1}b}=a.
 $$
  \item
{\em The Set Theoretic Yang--Baxter Property}:
 $$
(S\times id)(id\times S)(S\times id)=(id\times S)(S\times
id)(id\times S)\colon X^3\to X^3.
 $$
The Yang--Baxter property implies the following equations among the
operations.
 $$
a^{bc_b}=a^{cb^c},\, c_{ba^b}=c_{ab_a}\,\hbox{ and }\,
{b_a}^{c_a}={b^c}_{a^c}.
 $$
  \item
{\em The Biquandle Property}: In terms of the operations introduced
earlier, this is the property that
 $$
a^a=a_a
 $$
for all $a\in X$.
 \end{enumerate}

In~\cite{FJK} and elsewhere there are slightly different
conventions. The operations are defined directly by the switch $S$.
Suppose $S(a,b)=(b\vee a,a\wedge b)$. Then the operations $\vee,
\wedge$ are defined in terms of the above as
 $$
b\vee a=b_{a^{-1}},\, a\wedge b=a^{b_{a^{-1}}}.
 $$
The advantage of the new operations is that the formul\ae\ involved
in the Yang--Baxter equations and the biquandle condition become
much simpler.

The map, $G\colon X^2\to X^2$, defined by $G(a,b)=(b_a,a^b)$ is
called the {\em sideways} map in~\cite{FJK}. The biquandle property
says that $G$ preserves the diagonal in $X^2$.

Knot and link invariants can be defined from the above axioms by a
mix\rq{}n\rq{}match choice. Conditions 1 and 2 are usually
satisfied.

A {\em biquandle} satisfies all the axioms. A {\em birack} satisfies
all the axioms except the biquandle condition. {\em Quandles} and
{\em racks} are the same except the binary operation $a_b=a$ is
trivial.

\begin{figure}
  \centering\includegraphics[width=10cm]{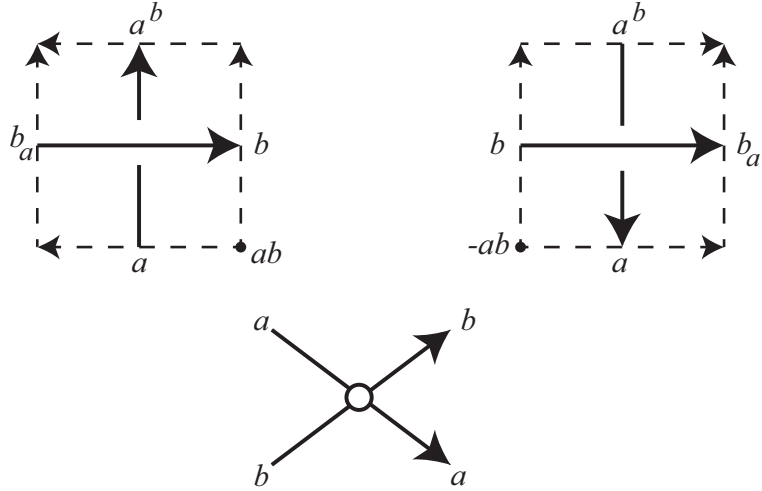}
  \caption{Conditions at a crossing}\label{fg:biqop}
 \end{figure}

So if all the edges of a virtual diagram are labelled by a biquandle
and satisfy the conditions at a crossing according to
Fig.~\ref{fg:biqop}. then we can obtain invariants of a virtual
knot. Note that we can think of a square around each crossing. The
label on the square is $\pm ab$ according to the labels on the edges
and the sign of the crossing.

Fig.~\ref{fg:biqinv} below give pictorial evidence of the invariance
under the Reidemeister moves. Full details can be obtained
in~\cite{Ftref}. Many examples of small size can be found
in~\cite{BF2}.

 \begin{figure}
  \centering\includegraphics[width=10cm]{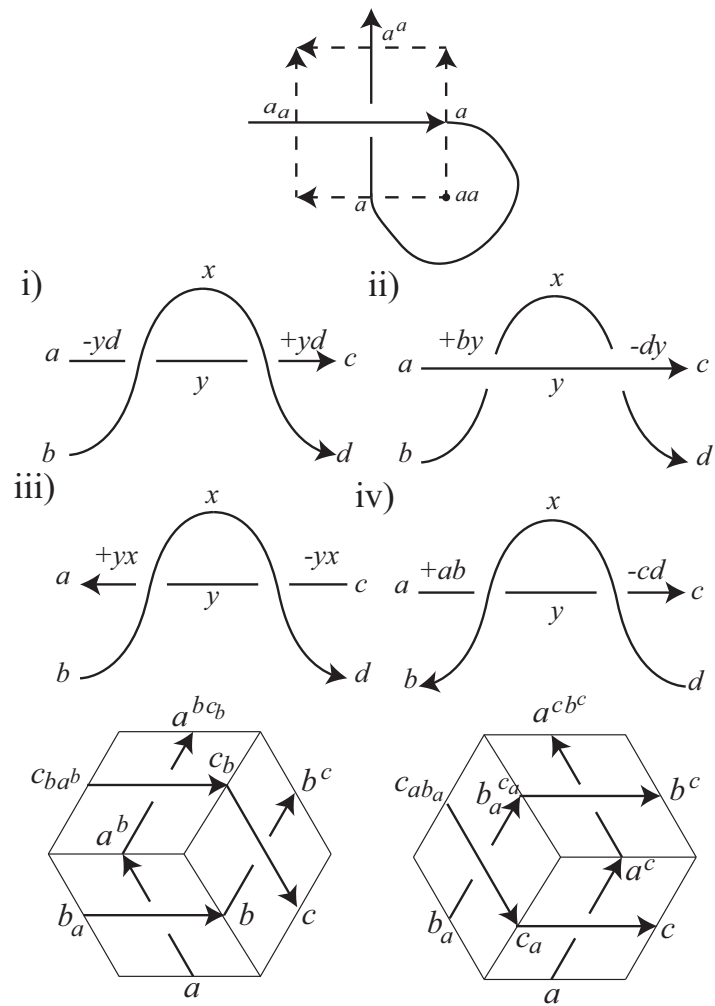}
  \caption{Invariance under the Reidemeister moves}\label{fg:biqinv}
 \end{figure}

On the other hand if we are looking for invariants of a doodle then
we would not expect the Yang--Baxter condition to be satisfied.

Suppose $R$ is an associative ring with a multiplicative identity,
1. If the set of labels is an $R$-module, then the birack,
biquandle, etc is called {\em linear} if there are are elements
$A,\,B,\,C,\,D$ of $R$ such that the switch map is defined by the
$2\times 2$ matrix
 $$
S=\left(\begin{array}{cc}A&B\\ C&D\end{array}\right).
 $$
There are many examples of linear biquandles with $R$ equal to the
quaternions, see for example~\cite{FBu,FBu2,GQ}.

{\em Affine} biquandles, generalizing linear ones need to be
investigated. One such is the Cheng labelling
 $$
S(x,y)=(y-1,x+1),
 $$
see~\cite{ZCheng}.

 \subsection{Homology of biracks and biquandles}

Let $X$ be the labels of a birack. A word $w=a_1\cdots a_n$ of
length $n$ is called an {\em $n$-cube}. For example a 1-cube is a
label, a 2-cube is a crossing and a 3-cube is a Reidemeister III
move.

For each $n$-cube there are $2n$ {\em faces} of dimension $n-1$,
 $$
\partial_i^- (a_1\cdots a_n)=a_1\cdots a_{i-1}a_{i+1}\cdots a_n
 $$
and
 $$
\partial_i^+(a_1\cdots a_n)=(a_1)^{a_i}\cdots (a_{i-1})^{a_i}(a_{i+1})_{a_i}\cdots (a_n)_{a_i}\,\hbox{for}\, i=1,\ldots,n
 $$
A cubical cell complex, $\Gamma X$, can be defined in the usual way
by identifying all the disjoint cubes having common faces.
Properties of this complex can now be used to define knot
invariants. For example the fundamental group has a presentation
 $$
\pi_1(\Gamma X)=\langle x,y\in X\mid  xy_x=yx^y\rangle.
 $$

Let $C_n=C_n(\Gamma X)$ be the free abelian group with basis the
$n$-cubes of $\Gamma X$. The homomorphism $\partial\colon C_n\to
C_n$ is defined  on cubes by
 $$
\partial=\sum_{i=1}^n(-1)^i(\partial_i^--\partial_i^+)
 $$
and extended linearly. The Yang--Baxter equations then imply that
the composition
 $$
C_n\to C_{n-1}\to C_{n-2}
 $$
is zero, and so we have a chain complex. The homology of this chain
complex $H_*(C)$ and the homology of any subchain complex is
therefore an invariant.

One important example is the degenerate subchain. A cube, $a_1\cdots
a_n$, is called {\em degenerate} if $a_i=a_{i+1}$, for some
$i=1,\ldots, n-1$. If the birack is a biquandle, so $a^a=a_a$ for
all $a$, then the degenerate cubes generate a subchain, $D_*$. The
quotient $BQ_*=C_*/D_* $ defines the biquandle chain complex and the
short exact sequence
 $$
0\to D_*\to C_*\to BQ_*\to 0
 $$
extends to a long exact homology sequence.

Another example is the {\em double} of a birack. Consider the set of
pairs of pairs $X^2\times X^2$. Then $W\subset X^2\times X^2$ is the
set of pairs of pairs defined by
 $$
W=\{ac, bc \mid a,b,c,\in X, c^a=c^b\}.
 $$
The doubled operations are
 $$
(ac)^{(bc)}=a^{b}c^b,\, (bc)_{(ac)}=b_ac^a\, \hbox{and}\,
G(ac,bc)=(b_ac^a,a^bc^a)
 $$

Doubling converts racks into biracks and quandles into biquandles.

The homology of the double of the 3-color quandle can be used to
distinguish the right and left trefoils, see~\cite{Ftref}.

Much work has been done on the homology of racks and quandles,
see~\cite{Carter,FRS2,FRS1,FRS4,FRS3}. However the homology of
biracks and biquandles is little understood and needs investigating.


\section{Graph-links}

\begin{figure}
  \centering\includegraphics[width=3.0in]{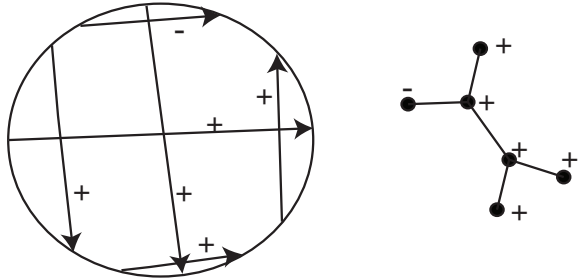}
  \caption{A Gauss diagram and its labeled intersection graph} \label{gaussdiagram}
 \end{figure}

It turns out that some information about the knot can be obtained
from a more combinatorial data: the {\em intersection graph} of a
Gauss diagram. The intersection graph is a graph without loops and
multiple edges, whose vertices are in one-to-one correspondence with
chords of the Gauss diagram. Two vertices of the intersection graph
are {\em adjacent} whenever the corresponding chords of the Gauss
diagram are {\em linked} (or intersect each other when they are
placed inside ), see Fig.~\ref{gaussdiagram}. Each vertex of the
intersection graph is endowed with the local writhe number of the
corresponding crossing.

\begin{figure}
\centering\includegraphics[width=4in]{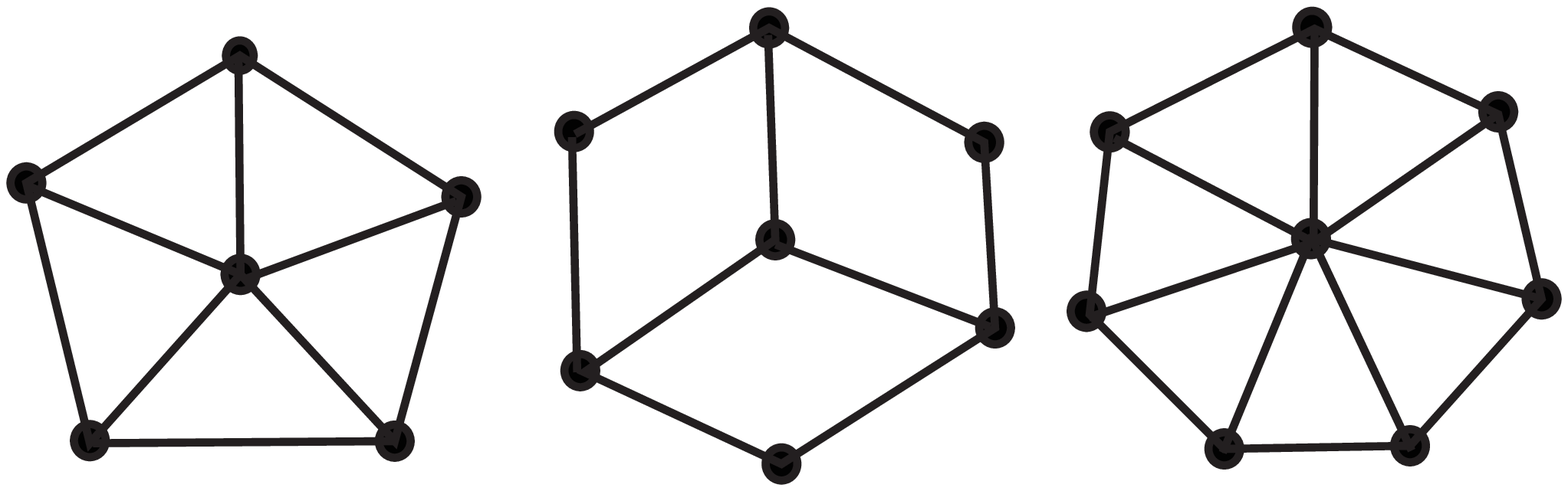}
\caption{Non-realizable Bouchet graphs} \label{Bouchet}
\end{figure}

However, sometimes a chord diagram can be obtained from the
intersection graph in a non-unique way, and some graphs (shown in
Fig.~\ref{Bouchet}) cannot be represented by chord diagrams at all.

Likewise virtual knots appear out of non-realizable chord diagram
and thus generalize classical knots (which have realizable chord
diagrams), graphs-links come out of intersection graphs: We may
consider graphs which are realized by chord diagrams, and, in turn,
by  virtual links, and pass to arbitrary simple graphs which
correspond to some mysterious objects generalizing links and virtual
links.

Traldi and Zulli~\cite{TZ} constructed a self-contained theory of
``non-realizable knots'' (the theory of {\em looped interlacement
graphs}) possessing lots of interesting knot theoretic properties by
using Gauss diagrams. These objects are equivalence classes of
(decorated) graphs modulo ``Reidemeister moves''.

\begin{figure}
\centering\includegraphics[width=4.0in]{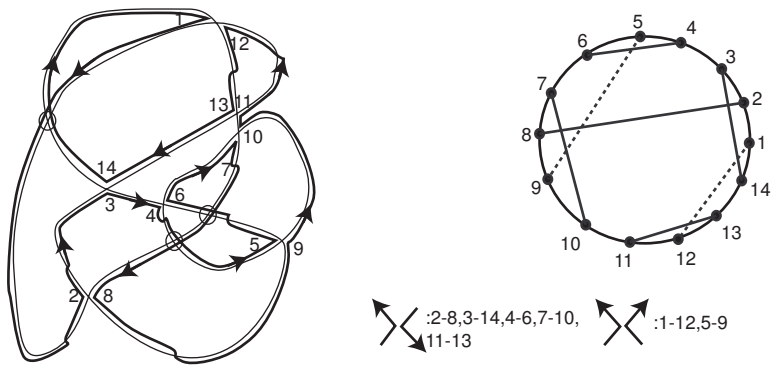} \caption{Rotating
circuit shown by a thick line; chord diagram} \label{rotations}
 \end{figure}

In~\cite{IM1} it was suggested another way of looking at knots and
links and generalizing them (the theory of {\em graph-links}):
whence a Gauss diagram corresponds to a transverse passage along a
knot, one may consider a {\em rotating circuit} which never goes
straight and always turns right or left at a classical crossing. One
can also encode the type of smoothing (Kauffman's $A$-smoothing or
Kauffman's $B$-smoothing) corresponding to the crossing where the
circuit turns right or left and never goes straight, see
Fig.~\ref{rotations}. We note that chords of diagrams are naturally
split into two sets: those corresponding to crossings where two
opposite directions correspond to {\em emanating} edges with respect
to the circuit and the other two correspond to {\em incoming} edges,
and those where we have two consecutive (opposite) edges one of
which is incoming and the other one is emanating.


After the two theories were constructed, some questions arose. It
was shown that there are graphs not being Reidemeister equivalent
(each theory has own Reidemeister moves) to an intersection graph of
a virtual knot diagram~\cite{IS1,Mant32}. The equivalence of the two
theories (the theory of looped interlacement graphs and the theory
of graph-knots) was proved in~\cite{Ily2}. Also, some invariants
were constructed, see~\cite{IM1,IM2,IM4,IS1,Nik1,TZ}.


\section{A List of Problems}


 \subsection{Problems in and related with virtual knot theory}

Below, we present a list of research problems closely connected with
virtual knot theory.

 \begin{enumerate}
  \item
{\em Recognising the Kishino Knot}: There have been invented many
ways to recognize the Kishino virtual knot (from the unknot): The
$3$--strand Jones polynomial, i.e.\ the Jones polynomial of the
$3$--strand cabling of the knot~\cite{KiSa}, the $\Xi$--polynomial,
see~\cite{Mant12}, the quaternionic biquandle~\cite{FB}, and the
surface bracket polynomial (Dye and Kauffman~\cite{DKMin}).
In~\cite{KADOKAMI} Kadokami proves the knot is non-trivial by
examining the immersion class of a shadow curve in genus two. One
can use the Manturov parity bracket~\cite{Mant31,Mant32,Mant35} to
show that essentially the flat Kishino diagram is its own invariant,
exhibiting the non-triviality of this knot. See~\cite{LKIntro} for
an exposition of this proof.

Are we done with this knot? Perhaps not. Other proofs of its
non-triviality may be illuminating. The fact that the Kishino
diagram is non-trivial and yet a connected sum of trivial virtual
knots suggests the question: {\em Classify when a non-trivial
virtual knot can be the connected sum of two trivial virtual knots.}
A key point here is that the connected sum of closed virtual knots
is not well defined and hence the different choices give some
interesting effects. With long virtual knots, the connected sum is
ordered but well-defined, and the last question is closely related
to the question of classifying the different long virtual knots
whose closures are equivalent to the unknot.
  \item
{\em Flat Virtuals}: Flat virtual knots, also known as {\em virtual
strings}~\cite{HR,TURAEV}, are difficult to classify. Find new
combinatorial invariants of flat virtual knots.

We would like to know more about the flat biquandle algebra. This
algebra is isomorphic to the Weyl algebra~\cite{FTur} and has no
(non-trivial) finite dimensional representations. An example that
goes beyond the usual restrictions is the (very simple) affine
biquandle used in \cite{Affine} to construct the Affine Index
polynomial invariant of virtual knots. One can make small examples
of the flat biquandle algebra, that detect some flat linking beyond
mod $2$ linking numbers, but the absence of other finite dimensional
representations presents a problem.
  \item
{\em The Flat Hierarchy}: The flat hierarchy is constructed for any
ordinal $\alpha$. We label flat crossings with members of this
ordinal. In a flat third Reidemeister move, a line with two $a$
labels can slide across a crossing labeled $b$ only if $a$ is
greater than $b.$ This generalizes the usual to theory of flat
virtual diagrams to a system with arbitrarily many different types
of flat crossings. Classify the diagrams in this hierarchy. This
concept is due to Kauffman (unpublished). A first step in working
with the flat hierarchy can be found in~\cite{Mant24}.
  \item
{\em Virtuals and the Theory of Doodles}: Compare flat theories of
virtual knots with theories of doodles. A doodle is represented by a
flat diagram in the plane. Reidemeister moves of type I and II are
allowed but not type III. So a triple point must not be allowed.
Khovanov has associated a group to doodles~\cite{KHO}. Commutator
identities can also be associated to a doodle. Cobordism of doodles
is defined by an immersed surface without triple points. Cobordism
classes represent elements of a free abelian group. There is
undoubtably a rich seam of results which could be found by
investigating doodles. For example there is no reason not to have
virtual crossings as well as flat crossings. This would give doodles
on a surface of higher genus, see~\cite{F,FTay,KHO}.
  \item
{\em Virtual Three Manifolds}: There is a theory of virtual
$3$--manifolds constructed as formal equivalence classes of virtual
diagrams modulo generalized Kirby moves, see~\cite{DK}. From this
point of view, there are two equivalences for ordinary
$3$--manifolds: homeomorphisms and virtual equivalence. Do these
equivalences coincide? That is, given two ordinary three manifolds,
presented by surgery on framed links $K$ and $L$, suppose that $K$
and $L$ are equivalent through the virtual Kirby calculus. Does this
imply that they are equivalent through the classical Kirby calculus?

{\em What is a virtual $3$-manifold?}: That is, give an
interpretation of these equivalence classes in the domain of
geometric topology. \bigbreak

Construct {\em another} theory of virtual $3$--manifolds by
performing surgery on  links in thickened surfaces $S_{g}\times
\mathbb{R}$ considered up to stabilization. Will this theory
coincide to that proposed by Kauffman and Dye~\cite{DK}?
  \item
{\em Welded Knots}: We would like to understand welded
knots~\cite{Rour,SATOH}. It is well known~\cite{KANENOBU,NELSON}
that if we admit forbidden moves to the virtual link diagrams, each
virtual knot can be transformed to the unknot. If we allow only one
forbidden move (e.g.\ the upper one), then there are lots of
different equivalence classes of knots. In fact the fundamental
group and the quandle of the virtual diagram are invariant under the
upper forbidden move. The resulting equivalence classes are called
{\em welded knots}. Similarly, {\em welded braids} were studied
in~\cite{FRR}, and every welded knot is the closure of a welded
braid. The question is to construct good invariants of welded knots
and, if possible, to classify them. In~\cite{SATOH} a mapping is
constructed from welded knots to ambient isotopy classes of
embeddings of tori (ribbon tori to be exact) in four dimensional
space, and it is proved that this mapping is an isomorphism from the
combinatorial fundamental group (in fact the quandle) of the welded
knot to the fundamental group of the complement of the corresponding
torus embedding in four-space. Is this Satoh mapping faithful from
equivalence classes of welded knots (links) to ambient isotopy
classes of ribbon torus embeddings in four-space?

Another interesting question is the following: If a welded knot has
trivial fundamental group, does this imply that the knot is trivial
as a welded knot?
  \item
{\em Long Knots and Long Flat Knots}: Enlarge the long knot
invariant structure proposed in~\cite{Mant17}. Can one get new
classical knot invariants from the approach in this paper? Bring
together the ideas from~\cite{Mant17} with the biquandle
construction from~\cite{GEN} to obtain more powerful invariants of
long knots. Long flat virtual knots can be studied via a powerful
remark due to Turaev (in conversation) to the effect that one can
associate to a given long flat virtual knot diagram $F$ a {\em
descending} diagram $D(F)$ (by always going over before going under
in resolving the flat (non-virtual) crossings in the diagram). The
long virtual knot type of $D(F)$ is an invariant of the long flat
knot $F.$ This means that one can apply any other invariant $I$ of
virtual knots that one likes to $D(F)$ and $I[D(F)]$ will be an
invariant of the long flat $F.$ It is quite interesting to do sample
calculations of such invariants and this situation underlines the
deeper problem of finding a full classification of long flat knots.

Bartholomew, Fenn, S.~Kamada and N.~Kamada have applied quaternion
invariants to long virtual knots, see~\cite{BFKK}.

See also Sec.~\ref{subsec:kar_chu} for more questions about long
virtual knots and long flat virtual knots.
  \item
{\em Virtual Biquandle}: Construct presentations of the virtual
biquandle with the a linear (non-commutative) representation at
classical crossings and some interesting structure at virtual
crossings. A start has been made by Bartholomew and Fenn~\cite{BF2}.
In this paper various biquandle decorations are made at classical
and virtual crossings which were found by a computer search.
 \item
{\em Virtual braids}: Is there a birack such that its action on
virtual braids is faithful?

Is the invariant of virtual braids in~\cite{Mant7}, see
also~\cite{Bardakov,KL,Mant19}, faithful?

The action defined by {\it linear} biquandles is not faithful. This
almost certainly means that the corresponding linear invariants of
virtual knots and links are not faithful~\cite{GQ} (see
also~\cite{BF2}).
  \item
{\em The Fundamental Biquandle}: Does the fundamental biquandle,
see~\cite{GEN} classify virtual links up to mirror images? (We know
that the biquandle has the same value on the orientation reversed
mirror image where the mirror stands perpendicular to the plane
(see~\cite{HR,HRK}).

Are there good examples of weak biquandles which are not strong?
This problem is solved in~\cite{Stan}.

We would like to know more about the algebra with 2 generators
$A,\,B$ and one relation $[B,(A-1)(A,B)]=0$ (see~\cite{FBu}). It is
associated to the linear case.
  \item
{\em Virtualization and Unit Jones Polynomial}: Suppose the knot $K$
is classical and not trivial. Suppose that ${\tilde K}$ (obtained
from $K$ by virtualizing a subset of its crossings) is not trivial
and has a unit Jones polynomial, $V(\tilde K)=1$. Is it possible
that ${\tilde K}$ is classical (i.e.\ isotopic through virtual
equivalence to a classical knot)?

Suppose $K$ is a virtual knot diagram with unit Jones polynomial. Is
$K$ equivalent to a classical diagram via virtual equivalence plus
crossing virtualization? (Recall that by crossing virtualization, we
mean flanking a classical crossing by two virtual crossings. This
operation does not affect the value of the Jones polynomial.)

Given two classical knots $K$ and $K',$ if $K$ can be obtained from
$K'$ by a combination of crossing virtualization and virtual
Reidemeister moves, then is $K$ classically equivalent to $K'?$

If the above two questions have affirmative answers, then the only
classical knot with unit Jones polynomial is the unknot.
  \item
{\em Virtual Quandle Homology}: Study virtual quandle homology in
analogy to quandle homology~\cite{CSH,GREENE}.
  \item
{\em Khovanov Homology}: Construct a generalization of the Khovanov
complex for the case of virtual knots that will work for arbitrary
virtual diagrams. Investigate the Khovanov homology constructed
in~\cite{Mant19,Mant21}. The main construction in this approach uses
an orientable atom condition to give a Khovanov homology over the
integers for large classes of virtual links. The import of our
question, is to investigate this structure and to possibly find a
way to do Khovanov homology for all virtual knots over the ring of
integers. Similar questions can be raised for the presently evolving
new classes of Khovanov homology theories related to other quantum
invariants (cf.~\cite{QK1,QK2}).

By a {\em K-full} virtual knot we mean a knot for which there exists
a diagram such that the leading (the lowest, or both) term comes
from the $B$-state. Analogously, one defines the {\em Kho-full} knot
relative to the Khovanov invariant. Call such diagrams optimal
diagrams. (It is easy to find knots which are neither K-full nor
Kho-full.)

Classify all K-full (Kho-full) knots.

Are optimal diagrams always minimal with respect to the number of
classical crossings?

Classify all diagram moves that preserve optimality.

Is it true that if a classical knot $K$ has minimal classical
diagram with $n$ crossings then any virtual diagram of $K$ has at
least $n$ classical crossings?

Can any virtual knot have torsion in the $B$-state of the Khovanov
homology (the genuine leading term of some diagram)? Here we use the
formulation of Khovanov homology given in~\cite{Mant19,Mant21}.

The behaviour of the lowest and the leading term of the Kauffman
bracket for virtual knots was studied in~\cite{Mant19}
and~\cite{Avdeev} and~\cite{NKamada}.
  \item
{\em Brauer algebra}: The appropriate domain for the virtual
recoupling theory is to place the Jones--Wenzl projectors in the
Brauer algebra. That is, when we add virtual crossings to the
Temperley--Lieb Algebra to obtain ``Virtual Temperley--Lieb
Algebra'' the result is the Brauer algebra of all connections from
$n$ points to $n$ points (see~\cite{TL,Fib2}). What is the structure
of the projectors in this context? Can a useful algebraic
generalization of the classical recoupling theory be formulated?
  \item
{\em Virtual Alternating Knots}: Define and classify alternating
virtual knots.

Find an analogue of the Tait flyping conjecture and prove it.
Compare~\cite{JZ}.

Classify all alternating weaves on surfaces (without stabilization).
 \item
{\em Crossing Number}: One of the most important problems in knot
theory is the problem of finding the minimal crossing number for a
given knot. It can be estimated by means of various invariants, the
simplest one, perhaps, being the Kauffman bracket estimate:
 $$
\mathrm{span}\langle K\rangle \leqslant 4n-4g,
 $$
Here $\mathrm{span}$ stands for the difference between the leading
degree and the lowest degree non-zero terms in the Laurent
polynomial defined by the Kauffman bracket, $n$ is the crossing
number, $g$ is the atom genus for the knot $K$ (also called  the
Turaev genus).

In the case $g=0$, we get the celebrated
Murasugi--Kauffman--Thistlethwaite's theorem on minimality of
alternating knot diagrams.

Is  there is a way to prove that the atom genus $g$ is minimal for a
concrete knot diagram? Then this together with a sharp estimate for
the span of the Kauffman bracket guarantees the minimality of the
diagram.

The estimate $4n-4g$ is exact for adequate diagrams (after
Thistlethwaite). It is exact for many diagrams which are not
adequate however: there are lots of minimal diagrams where the span
of the Kauffman bracket ``drops down'', for example, for torus
$(p,q)$ knots for $p>q>2$.

Why is the span of the Kauffman bracket smaller than expected for
such classes of knots?

One of possible explanations may be that the projection surface is
wrongly chosen. For a torus knot, it is more convenient to consider
it as being projected to the torus with no crossing points rather
than to the plane with many crossing points. Thus the following
problem arises: How do we use the Kauffman bracket, Thistlethwaite
spanning tree, atoms and other invariants to get estimates for the
number of crossings when projecting the knot diagram not to the
plane but rather to some surfaces of higher genera?

Certainly, when applying smoothing on the plane we get the final
expansion of the Kauffman bracket  as a linear combination of the
bracket for the unknot;  we shall get expression in terms of the
bracket of some ``basic'' torus knots, which makes the problem
harder.

Can we get estimates from Thistlethwaite's spanning tree and/or from
the Khovanov homology?

Which minimality theorems can be proved in this direction?

Is it possible to get an estimate for the atom genus when the atom
is considered not in the neighborhood of the plane but rather in a
neighborhood of the surface the knot is projected to and the
thickness of the Khovanov homology does not exceed $2+g$?
  \item
{\em Crossing number problems}: For each virtual link $L$, there are
three crossing numbers: the minimal number $C$ of classical
crossings, the minimal number $V$ of virtual crossings, and the
minimal total number $T$ of crossings for representatives of $L$.
There are also a number of unknotting numbers: The classical
unknotting number is the number of crossing switches needed to
unknot the knot (using any diagram for the knot). The {\em virtual
unknotting number} is the number of crossings one needs to convert
from classical to virtual (by direct flattening) in order to unknot
the knot (using any virtual diagram for the knot). Very little is
known. Find out more about the virtual unknotting number.

What is the relationship between the least number of virtual
crossings and the least genus in a surface representation of the
virtual knot.

Is it true that $T=V+L$?

Is there any algorithm for finding $V$ for some class of virtual
knots. For $T$, this is partially done for two classes of links:
quasialternating and some other, see~\cite{Mant9}. For classical
links and alternating diagrams see~\cite{MT,Mur}.

Are there some (non--trivial) upper and lower bounds for $T,\,V,\,L$
coming from virtual knot polynomials (see~\cite{LM})?
  \item
{\em Wild Virtuals}: Create the category of ``wild virtual knots''
and establish its axiomatics. In particular, one needs a theorem
that states when a wild equivalence of tame virtual links implies a
tame equivalence of these links.
  \item
{\em Vassiliev Invariants}: Understand the connection between
virtual knot polynomials and the Vassiliev knot invariants of
virtual knots (in Kauffman's sense). Some of that was done
in~\cite{GPV,VKT,Mant23,SAW2}.

The key question about this collection of invariants is this: {\em
Does every Vassiliev invariant of finite type, for classical knots
extend to an invariant of finite type for long virtual knots?} Here
we mean the problem in the sense of the formulation given
in~\cite{GPV}. In~\cite{VKT} it was pointed out that there is a
natural notion of Vassiliev invariants for virtual knots that has a
different notion of finite type from that given in~\cite{GPV}. This
alternate formulation needs further investigation.
  \item
{\em Embeddings of Surfaces}: Given a non-trivial virtual knot $K$.
Prove that there exists a minimal realization of $K$ in
$N=S_{g}\times I$ and an {\em unknotted} embedding of
$N\subset\mathbb{R}^{3}$ such that the obtained classical knot in
$\mathbb{R}^{3}$ is not trivial. (This problem is partially solved
by Dye in~\cite{D1}.
  \item
{\em Non-Commutativity and Long Knots}: It is known that any
classical long knot commutes with any long knot. This is definitely
not the case in the virtual category, see paper by Kamadas, Fenn et
al. However the commutativity of classical long knots may be
incorporated into the commutativity of virtual knots if the
following question is true. Is it true that if $K$ and $K'$ are long
knots and $K{\#} K'$ is isotopic to $K'{\#} K$ then there exists a
virtual long knot $L$, classical long knots $Q,\,Q'$, and
non-negative integer numbers $m,n$ such that
  $$
K=L^{m} {\#} Q, \quad K'=L^{n} {\#} Q',
  $$
where by $L^{m}$ we mean the connected sum of $m$ copies of the same
knot?
  \item\label{probl_old}
{\em The Rack Space}: The rack space was invented by Fenn, Rourke
and Sanderson~\cite{FENN,FRS2,FRS1,FRS3}. The homology of the rack
space has been considered by the above authors and Carter, Kamada,
Saito~\cite{CSH}. For low dimensions, the homology has the following
interesting interpretations. Two dimensional cycles are represented
by virtual link diagrams consistently colored by the rack, and three
dimensional cycles by the same but with the regions also colored.
See the thesis of Greene~\cite{GREENE}. So virtual links can give,
in this way, information about classical knots! For the second
homology of the dihedral rack, the results are given in Greene's
thesis. We now know that for a prime $p$ the third homology has a
factor $\mathbb{Z}_p$, see~\cite{NiebPrzy}.

Another line of enquiry is to look at properties of the birack
space~\cite{FRS2} and associated homology.

{\it Beyond problem number~{\em \ref{probl_old}} we list new
problems that are added to this revised problem list.}
  \item
Find new geometric/topological interpretations for the Jones
polynomial and for Khovanov homology.
  \item
Let $VKT/Z$ denote virtual knot theory modulo $Z$-equivalence, as
defined in the section above on virtual knot theory. Recall that two
virtual diagrams that are $Z$-equivalent have the same Jones
polynomial. It is known that classical knot theory embeds in virtual
knot theory. Does classical knot theory embed in $VKT/Z?$ That is,
suppose that $K$ and $K'$ are classical knot diagrams, and suppose
that $K$ and $K'$ are equivalent using virtual moves and $Z$-moves.
{\em Does it follow that $K$ and $K'$ are equivalent as classical
knots?}
  \item
Rotational virtual knot theory introduced in~\cite{VKT} is virtual
knot theory without the first virtual move (thus one does not allow
the addition or deletion of a virtual curl). This version of virtual
knot theory is significant because {\it all quantum link invariants
originally defined for classical links extend to rotational virtual
knot theory.} This theory has begun to be
explored~\cite{VKT,VKC,CVBraid} and deserves further exploration. We
have formulated a version of the bracket polynomial for rotational
virtuals that assigns variables according to the absolute value of
the Whitney degree of state curves. See Figure~\ref{rotate} for
examples of non-trivial rotational virtual links.

 \begin{figure}
  \centering\includegraphics[width=7cm]{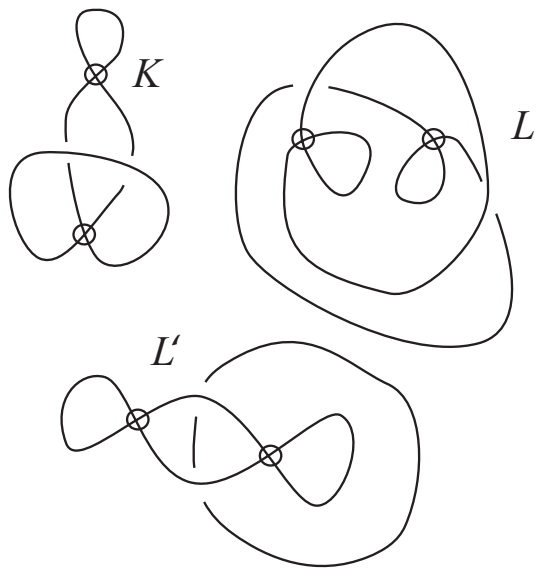}
  \caption{Examples of non-trivial rotational links}
  \label{rotate}
 \end{figure}

The rotational bracket can be generalised in the fashion
of the Arrow Polynomial. Beyond this there are all the quantum link
invariants and how they behave on rotational virtuals. We are in the
process of writing a new paper on rotational virtual knot theory.  A
natural class of invariants for rotational virtuals arises as
quantum invariants associated with finite dimensional
quasi-triangular Hopf algebras. This harks back to an early
project~\cite{KauRad,KRO,KRCAT,GEN} that needs further articulation.
In that work with Radford we articulate invariants that immediately
generalize to invariants of rotational virtual knots. These
invariants are defined via integrals on finite dimensional Hopf
algebras and should be studied for their own sake. We would also
like to know about the possibility of categorification of finite
dimensional Hopf algebras and, in particular, the meaning of
categorifying a right integral.
 \item
The arrow polynomial generalizes to a more powerful invariant of
virtual knots by defining a version of it for knots in specific
thickened surfaces.  The arrow polynomial for a knot in a thickened
surface has state variables that take into account the isotopy class
of the state curve and its arrow number as well. This gives a very
strong invariant of knots in thickened surfaces and it can be
applied to virtual knots by using a minimal surface representative
for the virtual knot. Investigating this {\em surface arrow
invariant} generalizes previous work on the surface bracket
polynomial~\cite{DKMin}.
 \item
Study concordance and cobordism invariants of virtual knots. In
particular, solve the question of virtual knots up to
pass-equivalence~\cite{FKT} (taking a direct generalization of
classical pass equivalence that gives the Arf invariant of knots).
Understanding this specific problem will promote our understanding
of concordance of virtual knots.
 \item
{\em Biquandles}: The first two problems are old chestnuts.

Give a descriptive representation of the free biquandle.

Give a topological explanation of the fundamental biquandle of a
knot.

Is there such a thing as a free partial biquandle: probably not.
  \item
Khovanov homology for virtual knots and links works directly with
mod-$2$ coefficients. With mod-2 coefficients there are no technical
difficulties associated with the fact that a virtual state curve can
resmooth to a single curve. Over the integers, this is a problem
that has to be dealt with.  An integral Khovanov homology for
virtual knots has been constructed by Manturov~\cite{Mant26}. It is
also the case that the work of Khovanov--Rozansky~\cite{KhR,KhoRoz}
categorifiying infinitely many specializations of the Homflypt
polynomial is an integral invariant for virtual knots. In principle,
the Khovanov--Rozansky work solves the problem of integral Khovanov
homology for virtuals, but to see this in terms of Khovanov's
original definition is a very good technical problem. Solve this
problem to clarify  the Manturov construction.
  \item
In~\cite{DKM,KaestKauff} we have studied a mod-2 categorification of
the arrow polynomial and discovered many baffling examples of pairs
of virtual knots that are discriminated by this link homology and
not distinguished by Khovanov homology mod-2 or the arrow
polynomial. We want to understand this new link homology and to that
end will do more computations and will try to understand the
structure of this homology theory for virtual links.
  \item
Make a systematic study of Vassiliev invariants for virtual knots
and links. There is ongoing work on this problem.
  \item
Generalize virtual knot theory to virtual 2-spheres in 4-space.
  \item
Find a way to effectively compute the Kauffman--Radford--Hennings
(KRH) invariants~\cite{KauRad} for three-manifolds via right
integrals on finite dimensional Hopf algebras~\cite{RI}. Find a way
to categorify the Kauffman--Radford--Hennings invariants for all
finite dimensional quasitriangular Hopf algebras. The KRH invariants
are actually invariants of rotational virtual links that are also,
after normalization, invariant under the moves of the Kirby
calculus. In this way, they give invariants of three manifolds and
virtual three manifolds. The formulation of the KRH invariants in
terms of integrals for finite dimensional Hopf algebras is very
elegant, but computation is very difficult. Thus these invariants
form a challenge for both virtual and classical knot theory. It is
worth pointing out again, that {\it the proper domain for all
quantum link invariants is rotational virtual knot theory.} Thus,
this problem about the KRH invariants is a test case for the
structure of quantum link invariants as a whole.
  \item
Create a combinatorial homotopy theory for Khovanov homology so that
the Khovanov homology of a knot or link is equivalent to the
homotopy type of an abstract complex (or category) associated with
the knot of link. This problem is meant in the sprit of Bar-Natan's
reformulation of Khovanov homology as an abstract chain homotopy
class of a categorical chain complex associated with the knot or
link. Recent work of Lifshitz and Sarkar~\cite{LSKho} goes very far
in this direction. Of course we would like a deeper connection, and
we would like to have the theory working for virtual knots and
links.

  \item
In~\cite{KNOTS} we give a formula for the Kauffman polynomial, due
to Jaeger, that expresses this invariant as a state sum over
oriented, partially smoothed links associated with the given link.
these states are evaluated by using the Homfly polynomial.
Generalize the Jaeger point of view to obtain a categorification of
the Kauffman polynomial by using chain complexes for these states
that are derived from the Khovanov--Rozansky categorification of the
Homflypt polynomial.
  \item
In~\cite{KNOTS,INT} we show how, by translating between knots and
planar graphs (the checkerboard and medial constructions) one can
associate a signed graph to a classical knot and that, with a choice
of nodes and an interpretation of the signs as generalized
electrical conductance, the conductance between two nodes of the
graph is an invariant of isotopy of the knot restricted to move in
the complement of these nodes. The moves on the graphs involve
replacements of pedant loops and edges, series and parallel
replacements and a star--triangle exchange move. These moves can be
applied to abstract (possibly non-planar) graphs. As a result we
obtain a generalization of knot theory (similar in spirit but quite
different from virtual knot theory) by taking the {\em electrical}
equivalence classes of signed graphs. Call this theory $EG_{\pm}.$
The conductance remains an invariant of these {\em Electric--Knots}
and shows that the theory is highly non-trivial. This formulation of
Electric Knots gives rise to many problems. Find new invariants of
Graph-Knots. Note that if we take a planar immersion of a signed
graph, then the induced cyclic orders of edges at the nodes of the
planar embedding give a natural way (analogous to the medial graph)
to associate a virtual knot diagram to the immersion. Explore this
relationship of Electric Knots and virtual knots.
  \item
It has been pointed out that the chain complex for Knot Floer
Homology~\cite{OS} (categorifying the Alexander--Conway polynomial)
is generated  by the states described in Formal Knot
Theory~\cite{FKT}. These states involve an assignment of pointers to
each region of the diagram (with two adjacent regions omitted) such
that each pointer marks a unique crossing in the diagram. The key
problem about this complex is that while the chain groups are
described via knot diagram combinatorics, the differential in the
complex seem to require high dimensional contact geometry. Even
though there is now a combinatorial definition of Knot Floer
homology via arc diagrams, this is an unsatisfactory state of
affairs. The problem is to find a combinatorial definition of the
differential on the complex generated by the Formal Knot Theory
states.
  \item
The Kauffman bracket polynomial of a virtual diagram can have the
leading term, i.e., the term having the highest possible degree,
equal to zero.

Assume that we have a diagram for which the Kauffman bracket
polynomial has the nonzero leading term. Classify moves
(compositions of the generalized Reidemeister moves) establishing
the equivalence in the class of diagrams with nonzero leading terms.
  \item
We still do not know whether the free knot whose Gauss diagram is a
heptagon (i.e.\ consists of 7 chords each of which is linked with
precisely two adjacent ones) is trivial.
  \item
Given a chord diagram $D$ we can construct the following formal
chain complex. Formally speaking, this complex will look like a
simplicial complex in all dimensions except zero, but one can get
rid of this by using duality arguments.

In dimension $-1$ there is just one chain denoted by $\emptyset$. In
dimension $0$, the chains are in one-to-one correspondence with
unoriented chords of the diagram $D$; in dimension $1$ the chains
are in one-to-one correspondence with pairs of unlinked chords, for
dimension $2$ we take triples of pairwise unlinked chords, etc.

The differential acts on the $k$-tuple, $\{a_{1},\dots, a_{k}\}$ of
pairwise unlinked chords to the $(k+1)$-tuple $\{a_{1},\dots,
a_{k},a_{l}\}$ obtained  by the addition of a chord unlinked with
all the previous ones. The sign is taken in a standard way to make
the complex well-defined with integral coefficients. Such a complex
is naturally defined by the intersection graph of the chord diagram.

 \begin{exa}
Let the chord diagram consist of $2n$ chords, where the the
following pairs of chords are linked: $(1,2)$, $(3,4)$, $(5,6),
\dots,(2n-1,2n)$. Then the corresponding simplicial complex is the
$(n-1)$-sphere (the $n$-th join of $S^{0}$.
 \end{exa}

Can we get any other simplicial complexes this way other than just
bouquets of spheres? An affirmative answer to this question may lead
to some homological calculations which are very easy for computer
implementation. This problem was motivated several years ago by an
attempt to construct the spectrification of the Khovanov homology.
It turns out that if some (virtual) knot in the  $A$-state has
exactly one circle then its Khovanov homology in the lowest quantum
grading looks exactly as described above. Thus, this approach may be
useful for understanding the Khovanov homology of chord diagrams  as
well as some other approaches to the Khovanov homology
spectrification obtained by recent work with Lipshitz.

On the other hand, this problem is motivated the construction of
graph-link theory. Indeed, every triangulated manifold admits a
shallow triangulation such that together with any 1-frame of a
$k$-simplex, it contains the whole simplex.

Considering the one-dimensional frame of such a triangulation as the
non-intersection graph of some chord diagrams, we get exactly the
chord diagram, whose lower term Khovanov homology looks as described
above. The problem is that not every graph is an intersection graph
(neither any graph is a non-intersection graph) of any chord
diagram. This leads to the graph-link theory for which all possible
manifolds appear in the lower degree homology theory. This leads to
some lower homotopy of graph links. The next problem can be
formulated as follows.

Given two diagrams of (virtual) knots or links with exactly one
circle in the $A$-state. Which combinations of Reidemeister moves on
the diagrams have to be be considered in such a way that whenever
two diagrams generate equivalent knots, there exists a chain of
moves from one diagram to the other with all intermediate  diagrams
having single-circle $A$-state and such that chord diagram homology
(and homotopy) corresponding to those intermediate diagrams do not
change under such moves?

It is worth comparing these thoughts with the recent work of
Lipschitz on the spectrification of Khovanov homology. Possibly, the
Bloom--Nikonov~\cite{Blo,Nik1,Nik2} approach to Khovanov
homology~\cite{ORS,OS} can be used here.
  \item
Which quantum invariants extend to virtual knots themselves without
restrictions? Certainly, there are such ones, e.g., the Kauffman
bracket and the Jones polynomial.
  \item
One can consider braids with even numbers of strands. Markov's moves
change the parity of the number of strands. Can one reformulate
Markov's theorem in such a way that only braids with even numbers of
strands take part in this formulation?
  \item
In his wonderful paper~\cite{Sar}, which firstly had the title
``Knot Floer Homotopy'', Sarkar constructs a cell complex. The
homology of this complex coincides with the Floer homology. In this
work the following main ingredients are used:
Manolescu--Ozsv\'{a}th--Sarkar--Szab\'{o}--Thurston approach
allowing one to define Heegaard--Floer homology by using grid
diagrams and shellability. Roughly speaking, if complexes look like
order complexes, then they can be replaced with cell (simplicial)
complexes homotopically equivalent to them.

Try to do the same things for Khovanov homology. In the initial
Khovanov definition this problem seems to be very difficult, since
the set of chains and states can hard be ordered. However, taking
into account Bloom's construction~\cite{Blo} this problem looks more
optimistic.
  \item
A complete invariant is used for proving that one theory is a part
of another theory. For example, the fact that the set of classical
knots is a part of the set of virtual knots, is true for the same
reason as a quandle with peripheral structure can be extended from
classical knots into virtual knots. Analogous statements can be done
for classical and virtual braids (the extension of Hurwitz's
action).

So, we get an actual problem: how can one construct a complete
invariant for virtual knots?
  \item
To prove that the invariant $\mathcal{F}$ constructed by
V.\,O.~Manturov for virtual braids, is complete for virtual braids
with more than two strands.
  \item
In~\cite{ManIly} it was showed how one could determine
non-invertibility of long virtual knots and non-commutativity of
long virtual knots. This method used two type crossings: an early
undercrossing and an early overcrossing. Further, for two types we
have different operations, these operations are similar but not
coincide (these operations are related to quandles). Regretfully,
this method gives nearly nothing for classical knots. But it is
universal: instead of the construction of quandle with two
operations one can consider other ways of constructing knot
invariants. In these ways we partition the set of crossings into two
subsets, and for each subset we use an operation. Moreover, two
operations are different enough to determine non-invertibility, and
in the same time are similar to give an invariant structure. In the
level of the Kauffman bracket polynomial this construction does not
work. We get two problems.

 \begin{enumerate}
  \item
Try to do the same for other objects than quandles, for example:
nonlinear quandles, quandles with rings having non-unique
decomposition of multiples.
  \item
Try to apply it on a ``categorified level'': use the usual chain
space for Khovanov homology and in it two different differentials
but commute with each other (these differentials correspond to
circle and star operations). These two different differentials could
arise from usual and odd Khovanov homologies in the case when we
have an identification (but maybe not canonical) between chain
spaces of these complexes. It is expected that these invariants
could recognize non-invertibility of knots.
 \end{enumerate}
   \item
It is well known that flat virtual knots are easily algorithmically
recognizable, see~\cite{HS}. The absence of a geometric approach to
the definition of free knots leads to the lack of the methods for
solving the following three problems.
   \begin{enumerate}
    \item
Is it true that (any) connected sum of free knots is trivial?
    \item
Are free knots algorithmically recognizable?
    \item
Prove that free knots do not commute in the general case.
   \end{enumerate}
In the first case we suggest that the answer is affirmative, and the
proof should be a certain modification of the analogous proof for
virtual knots (see, e.g.,~\cite{KM})

In the second case we suggest that the answer is negative. It seems
to us that free knots are a very complicated object, where one can
construct models for many logical constructions.

The third problem seems to be solvable by using (and, possibly,
extending) standard methods in the free knot theory.
   \item
The functorial mapping was constructed, by means of it we
constructed the map from the set of virtual knots to the set of
virtual knots with orientable atoms.

Can one construct an analogous map which projects the set of virtual
knots to the set of classical knots? For example, one can try to
find an index by using characteristic classes. Here only one
condition which arises under the third Reidemeister move should
hold.

Note that in the case of homologically trivial knots we can consider
a two-branched covering like a projection.
   \item
Turaev constructed the map from the set of long flat knots to the
set of long virtual knots. Can one construct any map from the set of
long free knots to the set of long virtual knots? This question can
be extended and one can try to construct any map which ``enlarge''
the structure.
   \item
The problem about cobordisms in sections: Let us have a free knot
and its cobordism. Construct a parity on the given free knot, which
is defined by using only this cobordism and respect only moves
inside this cobordism.
   \item
Prove or disprove the conjecture about the non-uniqueness of minimal
representative of a free link, i.e.\ there exists a free link having
several minimal representatives.
  \item
If $X$ is the free rack, is $\Gamma X$ a cat(0) space? A positive
answer would imply that all its higher homotopy groups are trivial.
 \item
Is there any algorithm for recognition whether two graph-links are
equivalent or not? Our conjecture is ``no''. The idea behind that is
that graph-links are complicated enough, and possibly, they contain
sufficiently many degrees of freedom to include something like
Turing machines. All mathematics is roughly split into the
recognizable one (low-dimensional topology, hyperbolic groups,
decidability) and the non-recognizable one (topology of dimension
$4$ and higher, arbitrary finitely presentable groups, Turing
machines etc). The usual argument for undecidability for finitely
presented group allows one to construct some universal (semi)groups
which includes the apparatus of the Turing machines. We think that
graph-links are enough complicated to include similar things: they
contain arbitrary graphs, and Reidemeister moves, in principle,
could play the role of ``rules for formal languages'' or ``group
relations''. It seems unlikely that Reidemeister moves collapse
graph-links to anything simple enough because of the parity
considerations: there are ``irreducibly odd'' graph-links which are
stable in the sense that every equivalent graph contains the initial
graph as a subgraph.
  \item
Can one construct a projection from the set of graph-links to the
set of realizable graph-links?
  \item
Construct a parity on graph-links by using Bouchet's criterion about
the realizability of a graph. Try to find a parity which is
responsible for the cyclically 6-edge connection, see~\cite{Bou1}.
  \item
Construct generalizations of the Frobenius
extension~\cite{Frobenius} and the Rasmussen for ``rigid''
graph-links with orientable atoms.
  \item
Is it true that two equivalent realizable graph-links are equivalent
in the class of realizable graph-links? If it is not true, then
construct an example.
  \item
Construct a group for graph-links which is analogous to the group
from~\cite[Subsection 8.9]{ManIly}.
  \item
Construct ``graph-braids''.
   \item
{\em Problems on Free Knot Cobordism}: The methods used for proving
the fact that the invariant $L$~\cite{IMN2} gives an obstruction to
the sliceness are not immediately generalized for obtaining lower
estimates on the slice genus of free knots. Two reasons are as
follows. First, when we define a {\em parity} and {\em justified
parity} for double lines on $\mathcal{D}$, we chose an arbitrary
curve connecting the two preimages of the point on the double curve.
We assert that any two curves connecting these two preimages (and
behaving correctly in neighborhoods of the ends) are homotopic. It
is true in the case of cobordism of genus zero, but in the case of a
surface of an arbitrary genus $h$ it is, indeed, not true.

Thus, in order to define a parity for double lines we have to impose
some restrictions on the spanning surface: we have to require that
the cohomology class dual to the graph $\Psi$ was
$\mathbb{Z}_{2}$--homologically trivial. The significance of the
property for even-valent graphs to be
$\mathbb{Z}_{2}$--homologically trivial is closely connected with
{\em atoms} (for details see Sec.~\ref{sec:atoms}
and~\cite{Mant32}).

Another problem is that the Reeb graph of an arbitrary Morse
function (not necessarily corresponding to the disc) is not
necessarily a disc. Consider Fig.~\ref{nottree}.

 \begin{figure}
  \centering\includegraphics[width=200pt]{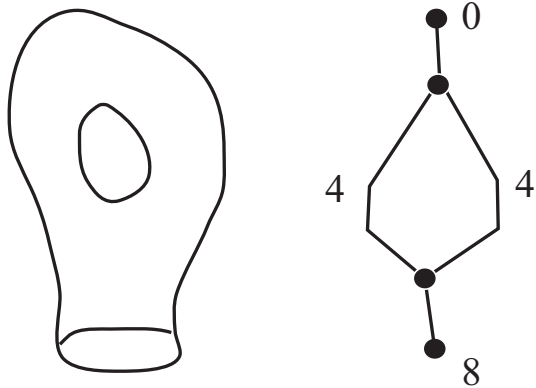}
  \caption{The Reeb graph corresponding a cobordism of genus one}\label{nottree}
 \end{figure}

Thus, starting from a free knot $K$ for which, we say, $L(K)=8$, we
(in principle) can turn it by a Morse bifurcation into free
two-component link consisting of two free knots  $K_{1}$ and
$K_{2}$, for which $L(K_{i})=4$, and then by another Morse
bifurcation we can reconstruct this trivial link into the unknot.
The invariant $L$ is not an obstruction to this, since the sum of
$4$ and $-4$ is zero.

In some cases we can overcome these two difficulties for cobordisms
(of arbitrary genus).

Let $\mathcal{D}_{g}$ be a surface with boundary $S^1$. Obviously,
the collection of double lines of $\mathcal{D}_{g}$ defines a
relative $\mathbb{Z}_{2}$--homology class $\kappa\in
H_1(\mathcal{D}_g,S^1;\mathbb{Z}_2)$. This homology class is an
obstruction for the surface to be checkerboard-colorable; also, this
is an obstruction for well-definedness of even/odd double lines.

Namely, if we look at the definition of an even/odd double line: we
see that there is an ambiguity in the choice of path connecting two
preimages of a generic point on the double line. For the case of a
disc cobordism, the parity of double lines is well defined, because
all such curves are homotopic. For $\mathcal{D}_{g}$ the unique
obstruction to this well-definedness is the class $\kappa$.

We call a cobordism of genus $g$ {\em checkerboard} (or {\em
atomic}\/) if the corresponding class $\kappa$ vanishes.

The next task (after detecting which $1$-stratum is even and which
one is odd) is to distinguish between $b$ and $b'$. To this end, one
should do the same for preimages of points lying on odd $1$-strata,
connect them by a generic curve, and count the intersection with
{\em even} double lines. So, we see that the only obstruction is the
relative $\mathbb{Z}_{2}$--homology class $\kappa'\in
H_1(\mathcal{D}_g,S^1;\mathbb{Z}_2)$ generated by even double lines.

We say that a checkerboard cobordism is {\em $2$-atomic} if
$\kappa'$ vanishes.

The following theorem holds.

 \begin{theorem}[see~\cite{Mant34}]
Assume for a $1$-component framed $4$-graph $K$ we have $L(K)\neq
0$. Then there is no $2$-atomic cobordism spanning the knot $K$ of
any genus.
 \end{theorem}

In the papers~\cite{ManMan,MM2} the fourth author constructed a
strengthening $G_{m}$ of the group $G$ given in the present work (in
the notation~\cite{ManMan} our group $G$ is $G_{1}$) and the
invariants of free knots with values in the classes of conjugate
elements from $G_{m}$. The idea is as follows. Even chords are
further partitioned into chords of different types, it leads to more
accurate definition of generators and relations in the group; these
constructions are closely connected with {\em iterated parities} and
the map deleting odd crossings. It seems that all invariants related
to the groups $G_{m}$ also give an obstruction to the sliceness of a
knot. Moreover, in the paper~\cite{Ma23} the author constructed
invariants of {\em virtual knots} in which the over/undercrossing
structure was taken into account besides the parity of chords. We
devote a separate paper to the investigation of a connection of
these invariants with cobordisms.
 \end{enumerate}

 \subsection{Problems on virtual knot cobordism}
The material in this section will appear in~\cite{KDK,VKC} and there
have been studies of virtual knot cobordism at the free knot level
by Manturov~\cite{Mant34,Mant35,Mant36}.

\begin{figure}
  \centering\includegraphics[width=10cm]{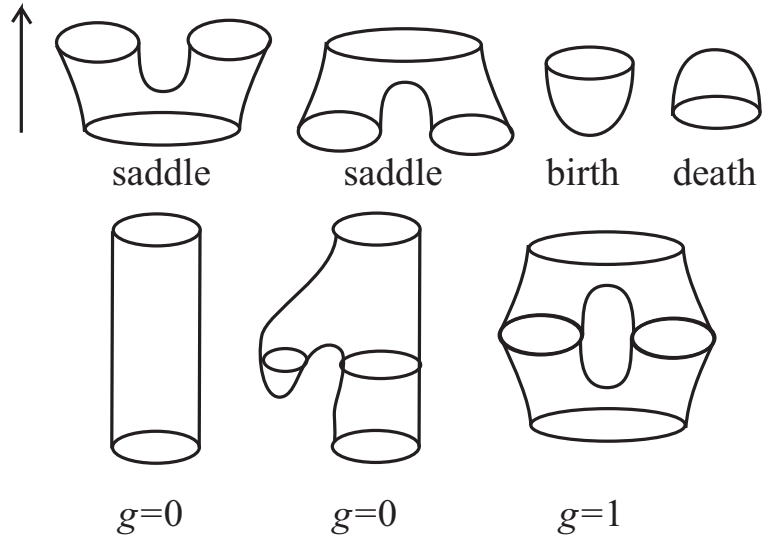}
  \caption{Saddles, births and deaths}\label{saddle}
 \end{figure}

 \begin{definition}
Two oriented knots or links $K$ and $K'$ are {\it virtually
cobordant}  if one may be obtained from the other by a sequence of
virtual isotopies (Reidemeister moves plus detour moves) plus
births, deaths and oriented saddle points, as illustrated in
Fig.~\ref{saddle}. A {\it birth} is the introduction into the
diagram of an isolated unknotted circle. A {\it death} is the
removal from the diagram of an isolated unknotted circle. A saddle
point move results from bringing oppositely oriented arcs into
proximity and resmoothing the resulting site to obtain two new
oppositely oriented arcs. See the figure for an illustration of the
process. Fig.~\ref{saddle} also illustrates the {\it schema} of
surfaces that are generated by  cobordism process. These are
abstract surfaces with well defined genus in terms of the sequence
of steps in the cobordism. In the figure we illustrate two examples
of genus zero, and one example of genus $1$. We say that a cobordism
has genus $g$ if its schema has that genus. Two knots are {\it
cocordant} if there is a cobordism of genus zero connecting them. A
virtual knot is said to be a {\it slice} knot if it is virtually
concordant to the unknot, or equivalently if it is virtually
concordant to the empty knot (The unknot is concordant to the empty
knot via one death.). This definition is based on the analogy of
virtual knot theory and classical knot theory. It should be compared
with~\cite{TurCob} in the category of virtual strings.
 \end{definition}

\begin{figure}
  \centering\includegraphics[width=10cm]{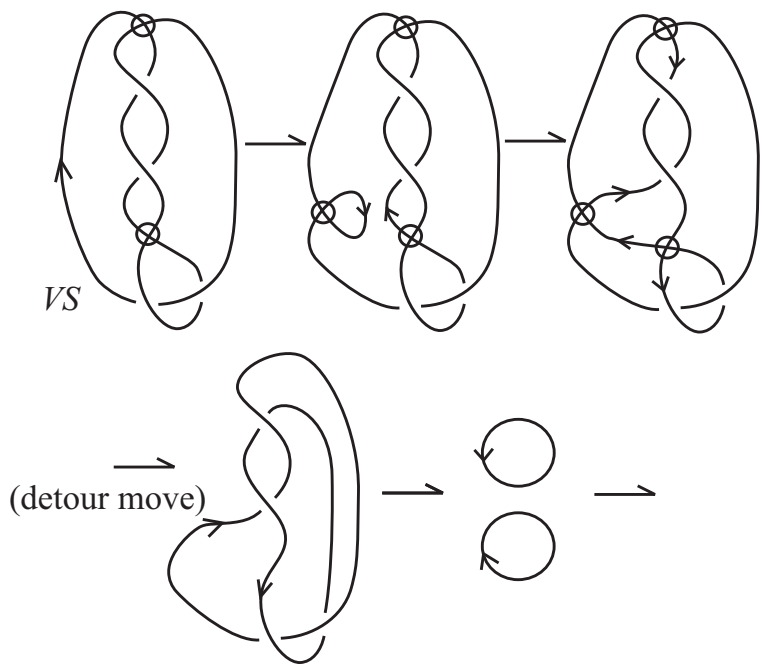}
  \caption{Virtual stevedore is slice}\label{vstevedore}
 \end{figure}

In Fig.~\ref{vstevedore} we illustrate the {\it virtual stevedore's
knot, VS} and show that it is a slice knot in the sense of the above
definition. We will use this example to illustrate our theory of
virtual knot cobordism, and the questions that we are investigating.
We define a virtual knot to be a {\it ribbon virtual knot} if it is
slice via only deaths and saddle point transformations. We do not
know at this point whether there are virtual slice knots that are
not ribbon. With the same definitions restricted to classical knots,
this is a long-standing open problem for classical knot theory.

\begin{figure}
  \centering\includegraphics[width=10cm]{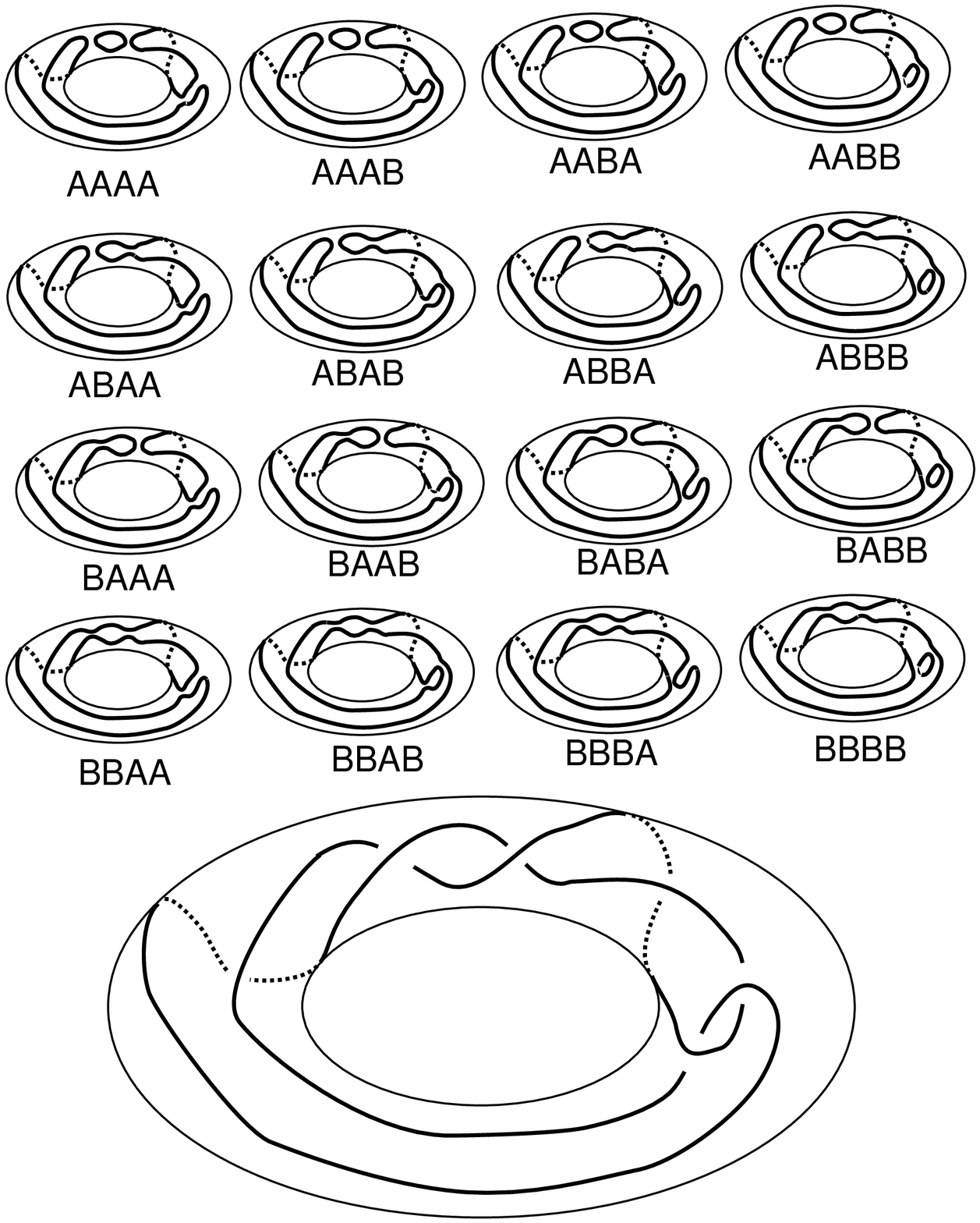}
  \caption{Virtual stevedore is not classical}\label{toral}
 \end{figure}

We first point out that the virtual stevedore ($VS$) is an example
that illustrates the viability of this theory. We can prove that
$VS$ is not classical by showing that it is represented on a surface
of genus one and no smaller. The technique for this is to use the
bracket expansion on a toral representative of $VS$ and examine the
structure of the state loops on that surface (see Fig.~\ref{toral}).
Note that in this figure the virtual crossings correspond to parts
of the diagram that loop around the torus, and do not weave on the
surface of the torus. An analysis of the homology classes of the
state loops shows that the knot cannot be isotoped off the handle
structure of the torus.

\begin{figure}
  \centering\includegraphics[width=10cm]{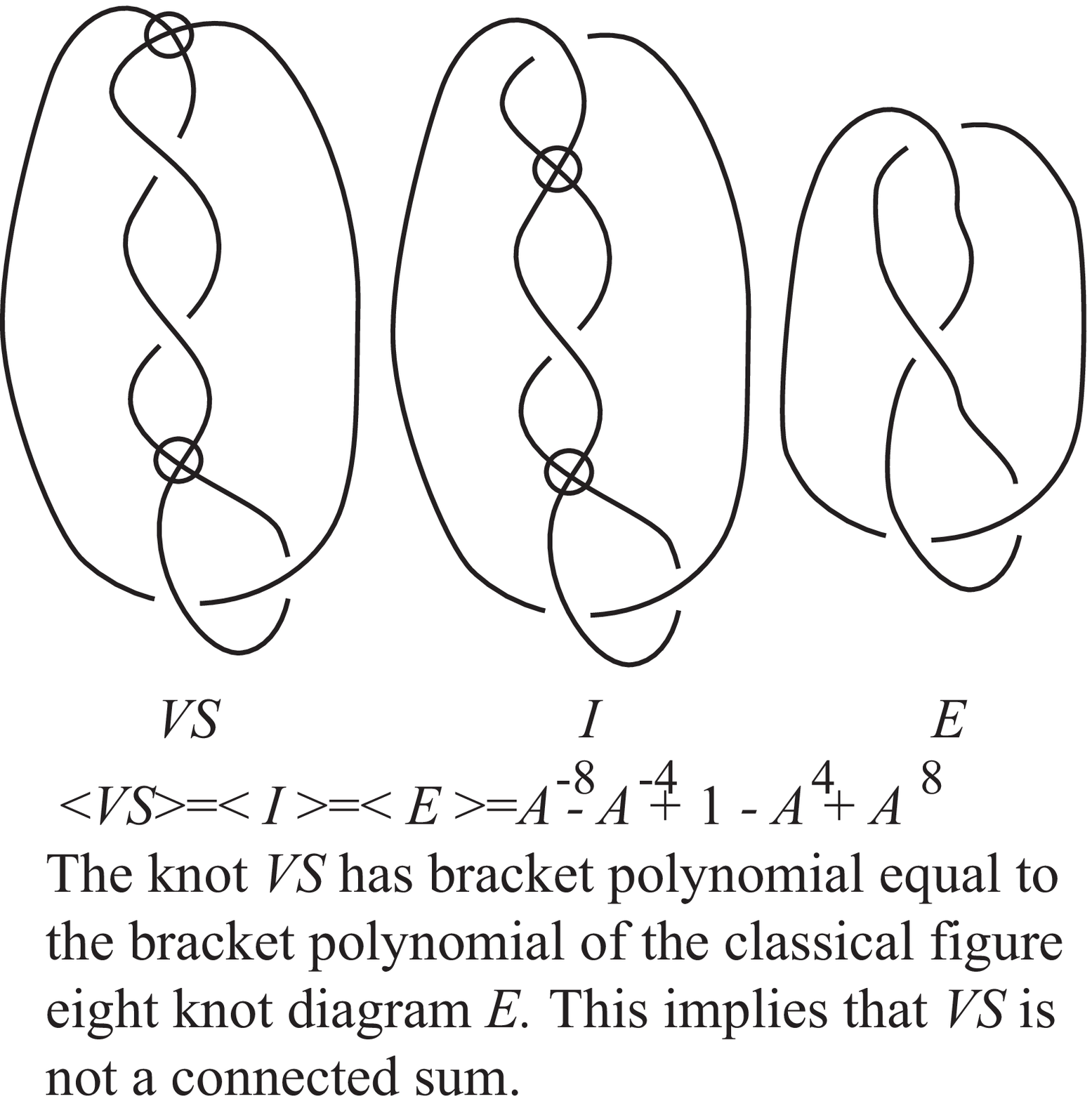}
  \caption{Bracket polynomial of the virtual stevedore}\label{VSbracket}
 \end{figure}

Next we examine the bracket polynomial of the virtual stevedore, and
show as in Fig.~\ref{VSbracket} that it has the same bracket
polynomial as the classical figure eight knot. The technique for
showing this is to use the basic bracket identity for a crossing
flanked by virtual crossings as discussed in the previous section.
This calculation shows that $VS$ is not a connected sum of two
virtual knots. Thus we know that $VS$ is a non-trivial example of a
virtual slice knot, {\it an example that constitutes proof in
principle that this project is viable.} We will face the problem of
the classification of virtual knots up to concordance.

 \subsubsection{Virtual surfaces in four-space}

\begin{figure}
  \centering\includegraphics[width=10cm]{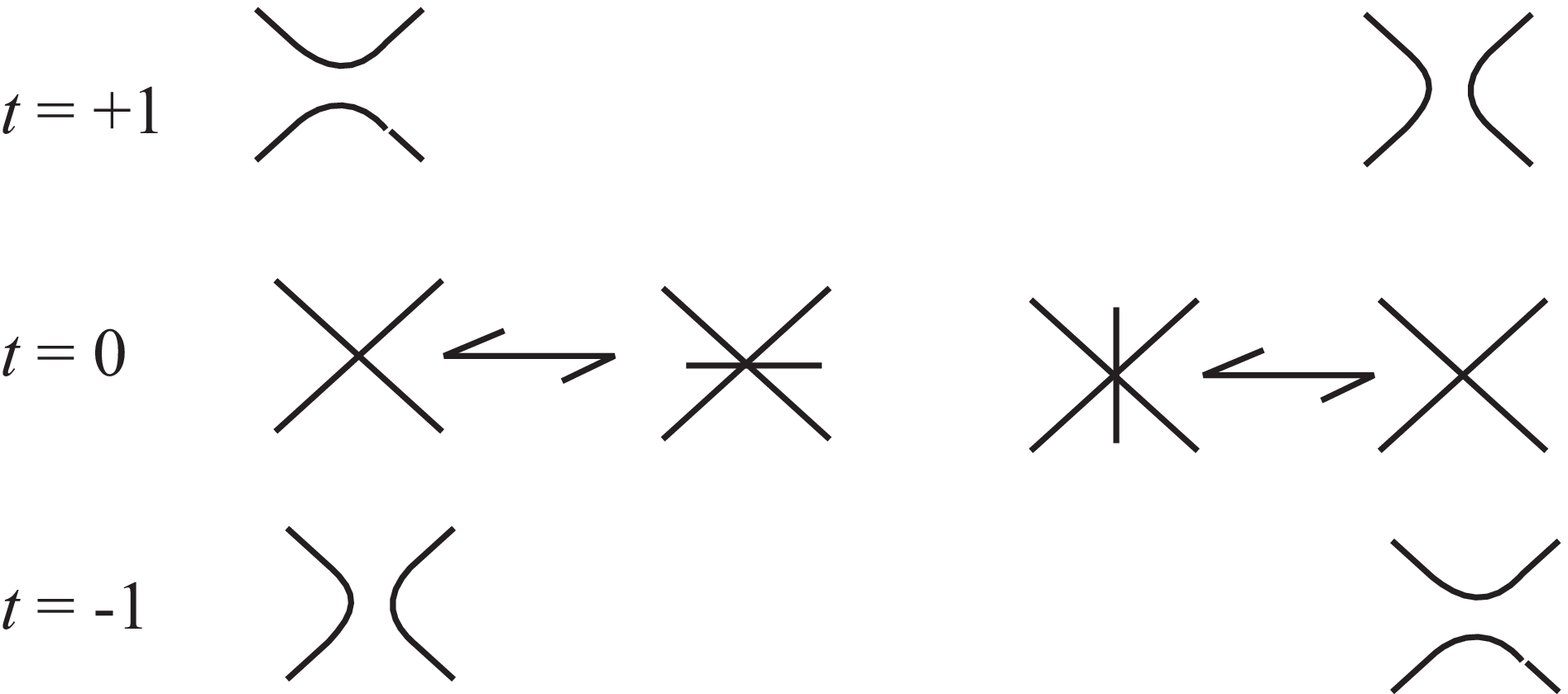}
  \caption{Middle level markers}\label{markers}
 \end{figure}

We now define a theory of virtual surfaces in four-space that is
given by moves on planar diagrams. One of the projects of this
proposal is to investigate the relationships between this
diagrammatic definition and more geometric approaches to virtual
$2$-knots due to Jonathan Schneider  and to Takeda~\cite{Takeda}
(see also~\cite{Swenton}). One should also compare with the
relationships to surfaces in four space for welded knots due to
Satoh~\cite{SATOH}. We make diagrammatic definitions as follows: We
use middle level markers as indicated in Fig.~\ref{markers} to
encode two directions of smoothing a marked crossing in a planar
diagram. The classical interpretation of such a marker is that it
represents a cobordism through a saddle point at the middle level
($t=0$ in the figure) where the forms of smoothing above ($t=1$) and
below ($t=-1$) are shown via the conventions in the figure. A
diagram with markers can then be interpreted as two cobordisms
attached at the middle. One cobordism goes downward to a collection
of possibly linked and knotted loops, the other goes upward to
another collection of linked and knotted loops. We will refer to
these as the {\it up-cobordism} and the {\it down-cobordism}. A
marked diagram is said to be {\it excellent} if both the up and the
down cobordisms end in collections of unlinked circles that can be
capped off with births (from the bottom) and deaths (at the top).
The resulting schema is then a two-sphere and classically represents
a two-sphere in four space. We take exactly this definition for a
{\it virtual two-sphere} where it is understood that the ends of the
two cobordisms will be trivial virtual links.

\begin{figure}
  \centering\includegraphics[width=10cm]{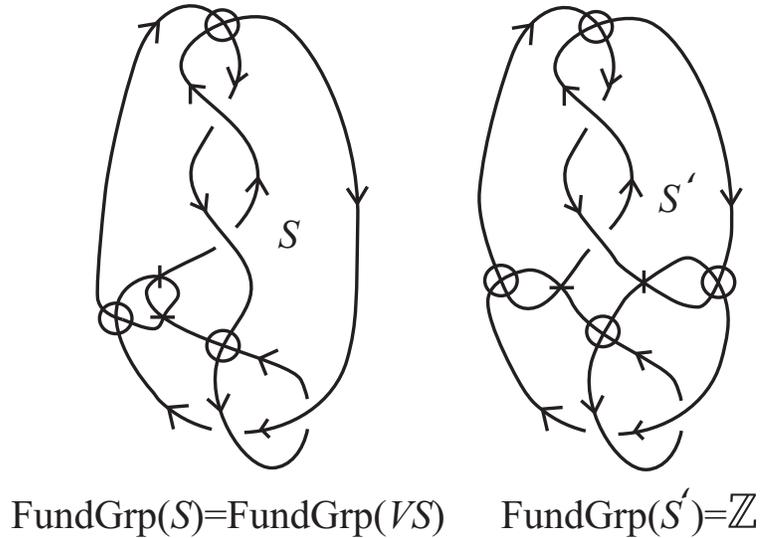}
  \caption{Two two-spheres}\label{twotwo}
 \end{figure}

Just as in classical theory, if a virtual knot is slice, then we can
make a virtual two-sphere from it by using the same cobordism both
up and down. Births for the original cobordism have to be
represented directly in the middle level. The slicing example for
$VS$, the virtual stevedore's knot, can be made into a two-sphere
this way. We show the middle level diagram for this sphere, called
$S$ in Fig.~\ref{twotwo}. In this same figure, we show another
middle level diagram for a virtual two-sphere $S'.$ In this case we
have used the fact (the reader can verify) that $VS$ can be sliced
from its right-hand side. The sphere $S'$ is obtained by slicing
upward from the left and downward from the right.

\begin{figure}
  \centering\includegraphics[width=10cm]{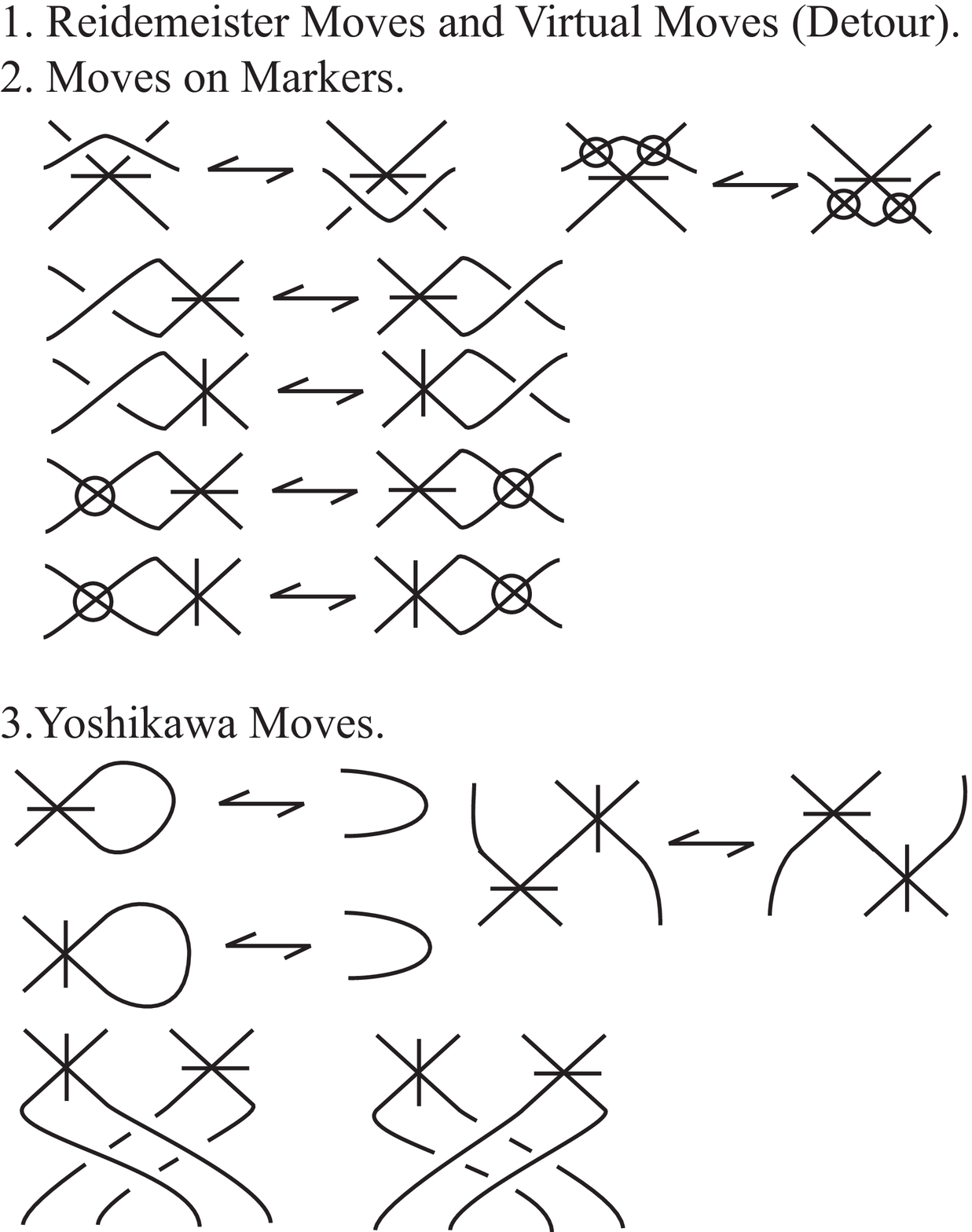}
  \caption{Middle level moves}\label{middlemoves}
 \end{figure}

\begin{figure}
  \centering\includegraphics[width=10cm]{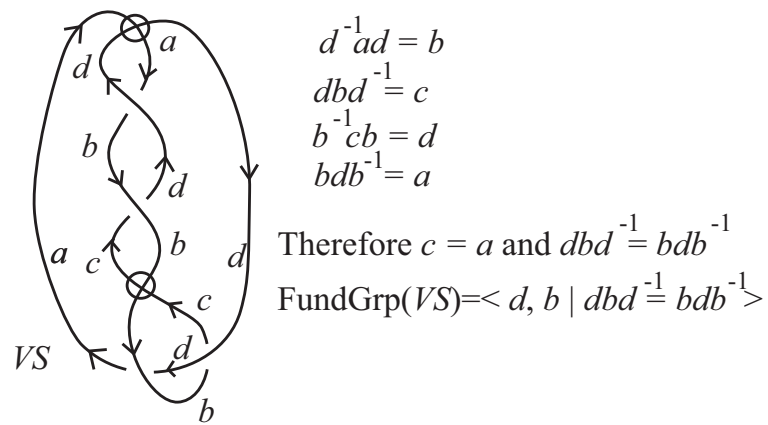}
  \caption{Fundamental group of $VS$}\label{fundgrp}
 \end{figure}

We give moves on the middle level diagrams to define isotopy of the
virtual two-spheres obtained from the middle level diagrams. The
moves are indicated in Fig.~\ref{middlemoves}. They are a virtual
generalization of the Yoshikawa moves that have been
studied~\cite{SYLee,Swenton}  for isotopies of the classical middle
level formulations. Thus we say the two two-spheres are {\it
isotopic} if one can be obtained from the other via these {\it
Generalized Yoshikawa moves}. In particular, the fundamental group
of the two-sphere,defined by adding relations at saddle points
exactly as in the classical case (but from the virtual knot
theoretic fundamental group) is an isotopy invariant. For example,
in Fig.~\ref{fundgrp} we calculate the fundamental group of $VS$ and
find that, in it the arcs whose elements must be identified to
obtain the fundamental group of the sphere $S$ of Fig.~\ref{twotwo}
are already identified in the fundamental group of $VS$. Thus we
find that the sphere $S$ is knotted since it has the same
non-trivial fundamental group as $VS.$ On the other hand, it is not
hard to see that the fundamental group of the sphere $S'$ is
isomorphic to the integers. At this writing we do not know if this
sphere is virtually unknotted.

The generalized Yoshikawa moves present a useful first formulation
for a theory of virtual surfaces. One of the advantages of this
approach is that we can adapt the generalization of the bracket
polyomial of Sang Youl Lee~\cite{SYLee} to obtain a bracket
invariant for virtual two-spheres. This will be an important subject
of investigation for this proposal. We want to know how this
diagrammatic formulation is related to immersions of surfaces in
four space that could represent virtual two-knots. In this case the
levels (movie of a cobordism) description that we have adopted gives
such an immersion, and one can begin the investigation at that
point. For these reasons, we believe that this formulation of
virtual cobordism and virtual surfaces will be very fruitful and
lead to many new results.

 \subsubsection{Virtual Khovanov homology}

Khovanov homology~\cite{Kh} (see
also~\cite{DB,BN2,KhoStat,KhoKauff,KM,Lee,Morrison,TuTu,Viro}) for
classical knots works, with mod-$2$ coefficients,  for  virtual
knots. It has been generalized by Manturov~\cite{Mant26,ManIly} for
a homology theory with integer coefficients. We have a new
formulation of Manturov's construction that simplifies some of the
choices in constructing the chain complex. We hope to see new
results from this technology. In particular we are examining the
structure of the Rasmussen invariant~\cite{Ras} with an eye to
generalizing it in this framework. Using the notions of cobordism
given here, we can define the {\it virtual four genus} of a virtual
knot as the least cobordism genus that it can attain. Thus slice
knots have genus zero and others will have higher genus in this
sense. The Rasmussen invariant gives a lower bound on the 4-ball
genus of classical knots. One can and we are investigating a similar
lower bound for the virtual four genus.

 \subsubsection{Band-passing and other problems}

The Arf invariant of a  classical knot can be interpreted as the
{\it pass-class} of the knot, where two knots are {\it pass
equivalent}~\cite{OK} if one can be obtained from the other by
ambient isotopy combined with switching pairs of oppositely oriented
pairs of parallel strands. The pass-class is a concordance invariant
of classical knots and closely related to the Alexander polynomial.
We would like to determine the pass-classes of virtual knots. This
problem appears difficult at this time due the lack of invariants of
the passing operation. We can obtain partial results by restricting
passing to only odd crossings (then the Manturov parity bracket
described in this proposal is an invariant of odd passing) but this
is only a step on the way to understanding the pass equivalence
relation for virtual knots. We expect that understanding this
relation will shed light on problems of knot concordance.

We cannot resist ending this collection of problems with a very
classical problem in skein theory. It is well-known that the
Alexander polynomial of a slice knot is of the form $f(t)f(1/t)$ up
to powers of $t$ and a sign. The problem is to give a proof of this
result using only Conway skein theory. This problem has resisted us
for a long time, and its resolution would surely shed light on
problems of knot cobordism classical and virtual.

 \subsubsection{Questions from Micah W.~Chrisman}

In recent talks, Kauffman has defined notions of cobordism,
concordance, slice, and ribbon for virtual knots.
Turaev~\cite{turaev_cob} has introduced a notion of cobordism for
knots in thickened surfaces. A knot in $\mathbb{R}^2 \times
\mathbb{R}$ that is slice in the classical is slice in both the
notions of Turaev and Kauffman. In~\cite{turaev_cob}, it is observed
that it is unknown whether the converse is true
  \begin{enumerate}
   \item
Are there any classical knots that are non-slice in the classical
sense but are slice in the the sense of Kauffman (when considered as
virtual knots) or Turaev (when considered as a knot in a thickened
surface)?
   \item
If a knot in the thickened surface is slice in the sense of Turaev,
does it necessarily stabilize to a virtual knot that is slice in the
sense of Kauffman?
  \end{enumerate}

 \subsection{Questions from V.~Bardakov}

The {\em virtual braid group} $VB_n$ is defined by the generators
$\sigma_1,\,\sigma_2,\ldots,\sigma_{n-1}$ which are generate the
braid group $B_n$ and  the generators $\rho_1,\,\rho_2, \ldots,
\rho_{n-1}$ which are generate the symmetric group $S_n$ also the
following mixed relations hold
$$
\sigma_i \rho_{i+1} \rho_i = \rho_{i+1} \rho_i \sigma_{i+1},~~~i =
1, 2, \ldots, n-2;
$$
$$
\sigma_i \rho_{j}  =  \rho_j \sigma_{i},~~~|i - j| > 1.
$$

 \begin{enumerate}
  \item
The group $G$ is called {\em linear} if there is an embedding of $G$
into the linear group $\mathrm{GL}_m(k)$ for some natural $m$ and a
field $k$.

Is it true that $VB_n$ is linear for all $n > 3$?

It is true for $n=3$ (unpublished result of V.~Bardakov and
P.~Bellingeri).
  \item
Let $F_{n+1} = \langle x_1, x_2, \ldots, x_n, y \rangle$ be a free
group of rank $n+1$. There is a representation $\rho\colon VB_n \to
\mathrm{Aut}(F_{n+1})$:
$$
\sigma_{i} : \left\{
\begin{array}{ll}
x_i \longmapsto x_i x_{i+1} x_i^{-1}, &  \\
x_{i+1} \longmapsto x_i, &  \\
x_k \longmapsto x_{k}, & \\
y \longmapsto y, &   \\
\end{array} \right.
~~~~ \rho_{i} : \left\{
\begin{array}{ll}
x_i \longmapsto y x_{i+1} y^{-1}, &  \\
x_{i+1} \longmapsto y^{-1} x_i y, &  \\
x_k \longmapsto x_{k}, &  \\
y \longmapsto y, &   \\
\end{array} \right.
$$
where $k \neq i,\,i+1$. This representation is an extension of the
Artin representation $\rho_A\colon B_n \to \mathrm{Aut}(F_{n})$
which is {\em faithful}, i.e.\ $\mathrm{ker}(\rho_A) = 1$.

Is the representation $\rho$ faithful for all $n > 2$?
  \item
The virtual knot theory is a subtheory of the {\em welded knot
theory}? The welded braid group $WB_n$ is the quotient of the
virtual braid group $VB_n$ by the relations
$$
 \rho_i \sigma_{i+1} \sigma_i = \sigma_{i+1} \sigma_i \rho_{i+1},~~~i = 1, 2, \ldots, n-2.
$$
The group $WB_n$ is a subgroup of $\mathrm{Aut}(F_n)$ and the
representation $\rho_A$ defined above is in fact the representation
of $VB_n$ into $WB_{n+1}$. In the virtual knot theory as in the
welded knot theory we have analogs of the Alexander theorem: every
link is equivalent to the closure of some braid.

Are there some non-equivalent virtual links $L_v$ and $L_v'$ that
are the closures of the virtual braids: $L_v = \widehat{\beta_v}$
and $L_v' = \widehat{\beta_v'}$ for some $\beta_v \in VB_n$,
$\beta_v' \in VB_m$ such that the welded links
$\widehat{\rho_A(\beta_v})$ $\widehat{\rho_A(\beta_v'})$ are
equivalent as welded links?
  \item
The {\em flat virtual braid group} $FVB_n$ is the quotient of the
virtual braid group $VB_n$ by the relations
$$
\sigma_i^2 = 1,~~~i = 1, 2, \ldots, n-1.
$$

Is there a representation $\varphi\colon FVB_n \to
\mathrm{Aut}(F_{m})$ for some $m$ such that in the image the
forbidden relations
$$
\varphi(\rho_i) \varphi(\sigma_{i+1}) \varphi(\sigma_i) =
\varphi(\sigma_{i+1}) \varphi(\sigma_i) \varphi(\rho_{i+1}),~~~i =
1, 2, \ldots, n-2
$$
do not hold?
 \end{enumerate}

 \subsection{Questions from Karene Chu}\label{subsec:kar_chu}

\begin{figure}
  \centering\includegraphics[height=5cm]{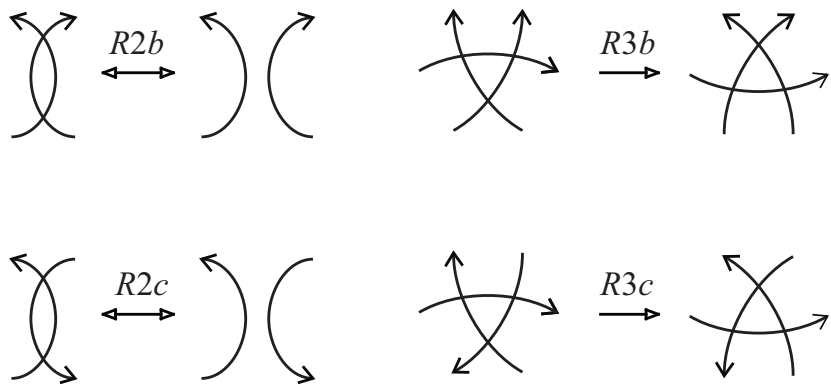}
  \caption{$R2b$, $R2c$, $R3b$ and $R3c$ moves}\label{fig:R2R3bc}
  \end{figure}

  \begin{enumerate}
   \item
{\em Long Virtual Knots and Long Flat Virtual Knots}\/: Even though
(round) flat virtual knots are not well understood, it turns out
long flat virtual knots can be completely classified. They are in
bijection with the subset of all signed permutations of in which
consecutive pairs are not sent to oppositely signed consecutive
pairs~\cite{Chu:FlatKnots}. For example, $1\mapsto (3,+) $,
$2\mapsto (2,-)$, $3\mapsto (4,-)$, $4\mapsto (1,+)$ is excluded
since the pair $(1,2)$ is mapped to the oppositely signed
consecutive pair $((3,+),  (2,-))$.  This classification is an
invariant on long virtual knots; in particular, it gives necessary
conditions for a long virtual which is classical, and can
distinguish the long Kishino knots from the unknot. How can we use
this invariant to understand long free knots, or round flat virtual
knots, both of which are quotients of long flat virtual knots?

Long flat virtual knots are also equivalent to long descending
virtual knots by sending each flat real crossing to a descending
crossing.  The long descending virtual knots all have cyclic
fundamental groups. What are the characteristics of the knot
polynomials on this subset of long virtual knots?

In fact, flat virtual pure tangles, whose skeletons are labelled
long strands, are also completely classified by the same approach.
There is a well-defined map from flat virtual braids in to flat
virtual pure tangles simply by interpreting the braid diagram as a
tangle diagram.  It was pointed out in Chu's talk that if this map
is injective, then this classification also gives normal forms for
flat virtual braids.  How to show that this map is injective?
 \item
{\em Virtual braids and flat virtual braids and braid-like virtual
knots}\/: The classifying spaces for both the flat virtual braid
group and the virtual braid group have been constructed
in~\cite{BEER:VirtualBraids}. (In this paper the flat virtual and
virtual braid groups are called the triangular and quasi-triangular
groups respectively.)  Indeed, the classifying space of the flat
virtual braid group is a beautiful space which is the quotient of
the permutohedron by a natural action of the symmetric group.  But
from another perspective, what are the characterizations of the flat
virtual or virtual braids in terms of embeddings into surfaces or
thickened surfaces with boundary modulo stabilization?

The equivalence generated by the subset of ``braid-like''
Reidemeister II and III moves, defined in Fig.~\ref{fig:R2R3bc}, is
a proper subset of the equivalence generated by the set of all
Reidemeister II and III moves. ``Braid-like'' or ``acyclic''
Reidemeister II and III moves, shown as $R2b$ and $R3b$, are moves
whose diagrams contain no oriented cycles (with the crossing as
vertices).  Cyclic moves, $R2c$ and $R3c$, contain oriented cycles.
The set of virtual knot diagrams modulo only braid-like Reidemeister
moves are called braid-like virtual knots by Bar-Natan.  One can
construct quantum invariants of these simply by ``attaching
R-matrices of quantum groups to the real crossings'' and ignoring
cups and caps, unlike for classical knots (see~\cite{KhoLau}). What
are braid-like virtual knots topologically?
 \item
{\em Welded braids and Knots}\/: The welded braid group is the
basis-conjugating subgroup of the automorphism group of the free
group and includes the classical braids.   It has been completely
classified with normal forms in~\cite{GK}. We can ask which
properties of the classical braid group are also possessed by the
welded braid group, e.g., automaticity, orderability, linearity?
 \end{enumerate}

 \subsection{Questions from Micah W.~Chrisman}

 \subsubsection{Finite-type invariants}

Finite-type invariants of classical knots can be understood through
the Kontsevich integral, configuration space integrals, and
combinatorial formulae (see~\cite{KFI}). There is no corresponding
general theory for Vassiliev invariants of virtual knots.

There are two notions of finite-type invariants of virtual knots:
the Goussarov--Polyak--Viro finite-type invariants and the Vassiliev
finite-type invariants (introduced to the virtual knot case by
Kauffman). There is a universal finite-type invariant for the GPV
finite-type invariants by way of the Polyak algebra. On the other
hand, there is no known universal finite-type invariant for
Vassiliev finite-type invariants.

In fact, the two types are quite different. The GPV finite-type
invariants are finitely-generated at every order. The Vassiliev
finite-type invariants are however infinitely generated. There is
even an infinite number of different extensions of the Conway
polynomial to virtual knots, all satisfying the same skein relation
and whose coefficients give an infinite family of finite-type
invariants of every order~\cite{MC_lattice}. Also, there are
obstructions for when a Vassiliev finite-type invariant can be
represented by an arrow diagram formula in the sense of
Goussarov--Polyak--Viro. For example, an integer valued Vassiliev
invariant which is invariant under the virtualization move cannot be
a GPV finite-type invariant~\cite{MC_virtmove}.
 \begin{enumerate}
  \item
Classify all obstructions for a Vassiliev finite-type invariant to
be a GPV finite-type invariant.
  \item (see also~\cite{CD_three})
Let $\Sigma$ be a compact oriented surface. Let $G$ be an abelian
group, $\mathscr{VK}$ the isotopy classes of oriented virtual knots,
and $\mathscr{K}[\Sigma]$ the isotopy classes of oriented knots in
$\Sigma \times I$. Let $\kappa\colon \mathscr{K}[\Sigma] \to
\mathscr{VK}$ denote the projection of knots on $\Sigma \times I$ to
virtual knots. For a set $S$, let $\mathbb{Z}[S]$ denote the free
abelian group generated by $S$. If $v\colon\mathscr{K} \to G$ is a
Vassiliev finite-type invariant of virtual knots, prove that there
is a finite-type invariant of knots in $\Sigma \times I$ such that
the following diagram commutes:
\[
\xymatrix{\mathbb{Z}[\mathscr{K}[\Sigma]] \ar@{-->}[dr]^{V}
\ar[d]_{\kappa} &\\ \mathbb{Z}[\mathscr{VK}] \ar[r]_{v} & G.}
\]
 \end{enumerate}

 \subsubsection{Prime decompositions of long virtual knots}

Long classical knots commute under the the operation of
concatenation. On the other hand, concatenation is not a commutative
operation for long virtual knots. This was first proved by Manturov
\cite{manturov_compact_long}, but has since been established by many
authors using a variety of techniques. Certainly, if either $A$ or
$B$ is classical, then $A\#B \leftrightharpoons B\# A$. If
$A,\,B,\,C,\,D$ are linear prime non-classical long virtual knots,
$A\#B \leftrightharpoons C\#D$, and $A \# B$ is non-classical, then
$A \leftrightharpoons C$ and $B \leftrightharpoons D$
\cite{MC_prime}. This gives a partial answer to the Manturov's
conjecture. In the virtual knot case, an existence and uniqueness
theorem for prime decompositions of virtual knots has been
proved~\cite{koramat,matveevprime}.
  \begin{enumerate}
   \item
Besides concatenation, there are other kinds of decompositions for
long virtual knots. For non-classical long virtual knots $A$ and $B$
that are prime relative to these decompositions,  is it possible for
the concatenation $A \# B$ to be expressed as any other type of
decomposition into $A$ and $B$?
   \item
Prove an existence and uniqueness theorem for prime decompositions
of long virtual knots.
   \item
Prove a generalization of Kuperberg's theorem to virtual tangles.
There are a few generalizations of Kuperberg's theorem appearing
explicitly in the
literature~\cite{twisted_links,MC_prime,KM1,roots}.
 \end{enumerate}

 \subsubsection{Connections between classical and virtual knot
theory}

Virtual knot theory can be applied to studying 1-knots in
3-manifolds and links in $S^3$ with at least two components. The
technique is called the {\em theory of virtual
covers}~\cite{MC,cm_fiber,km_fiber2}.

The technique may be described as follows. Let $N$ be a compact
oriented $3$-manifold and $K$ a smooth knot in $N$ (written $K^N$).
Suppose that $N$ admits a covering space $\Pi:\Sigma \times (0,1)
\to N$ and that there is a smooth knot $\mathfrak{k}$ in $\Sigma
\times (0,1)$ such that $\Pi(\mathfrak{k})=K$. Then $\mathfrak{k}$
can be considered as a knot in $\Sigma \times [0,1]$ via inclusion
and projected to a virtual knot $\upsilon$. Then we say that
$(\mathfrak{k}^{\Sigma \times (0,1)},\Pi,K^N)$ is a {\em virtual
cover with associated virtual knot} $\upsilon$. If $\Pi$ is regular
and $\mathfrak{k}$ is contained in a fundamental region of $\Pi$,
then it is said to be a {\em fundamental virtual cover}.

The theory has been studied in detail in the case that $N$ is the
complement of a fibered knot, a multi-component fibered link, and a
virtually fibered $3$-manifold (possesses a finite-index cover that
is fibered over $S^1$). For such $N$, inequivalent knots can be
detected by applying virtual knot invariants. The non-invertibility
of knots in manifolds can also be detected. Moreover, minimality
theorems can be proved for diagrams of knots in a manifold relative
to fiber of the fibration~\cite{km_fiber2}.
  \begin{enumerate}
   \item
The theory of virtual covers uses virtual knot theory to study links
of at least two components in $S^3$. Is there an extension of the
theory of virtual covers to one component links (i.e. knots) in
$S^3$.
   \item
Use virtual covers to find a relation between the Milnor
$\mu$-invariants and invariants of virtual knots. Can virtual covers
be used to apply Milnor's theory to one component links (i.e.
knots)?
   \item
How well does a virtual knot associated to a virtual cover of a
classical link behave under local moves (non-isotopy) applied to the
link? Of particular interest would be moves such as crossing
changes, self-$\Delta$ moves, and $Z$-moves.
   \item
Can virtual covers be used to detect mutations of classical links?
   \item
Every virtual knot $\upsilon$ is the associated virtual knot for a
fundamental virtual cover of some two component link $K \sqcup J$,
where $J$ is fibered and $\text{lk}(J,K)=0$. Find a complete set of
moves on such two component links which do not change the virtual
isotopy type of the associated virtual knot.
   \item
Let $J$ be a link in $S^3$ and $K$ a knot in the complement $N$ of
$J$. Let $\upsilon$ be a virtual knot associated to a virtual cover
$(\mathfrak{k}^{\Sigma \times (0,1)},\Pi,K^N)$. If $\upsilon$ is a
prime virtual knot, must $K$ be prime as a knot in $S^3$?
\end{enumerate}



\end{document}